\documentclass{elsarticle}
\IfFileExists{microtype.sty}{\usepackage{microtype}}{\relax}
\usepackage[british]{babel}
\usepackage[cmex10]{amsmath}
\usepackage{amsfonts}
\usepackage{amssymb}
\usepackage{amsthm}
\usepackage{enumerate}

\newcommand{\wt}{\mathrm{wt}\,}
\newcommand{\FF}{\mathbb F}

\newcommand{\KK}{\mathbb K}

\newcommand{\cG}{\mathcal G}

\newcommand{\GG}{\mathbb G}

\newcommand{\HH}{\mathbb H}
\newcommand{\cC}{\mathcal C}

\newcommand{\cH}{\mathcal H}
\newcommand{\cN}{\mathcal N}
\newcommand{\Rad}{\mathrm{Rad}\,}

\newcommand{\PG}{\mathrm{PG}}

\newcommand{\fA}{\mathfrak A}
\newcommand{\fB}{\mathfrak B}
\newcommand{\fC}{\mathfrak C}

\newcommand{\cod}{\operatorname{codim}}

\newcommand{\rank}{\mathrm{rank}\,}
\newcommand{\Fix}{\mathrm{Fix}\,}

\newtheorem{theorem}{Theorem}[section]
\newtheorem{lemma}[theorem]{Lemma}
\newtheorem{corollary}[theorem]{Corollary}
\newtheorem{prop}[theorem]{Proposition}
\theoremstyle{definition}
\newtheorem{remark}[theorem]{Remark}

\begin{document}
\begin{frontmatter}
\title{Line Hermitian Grassmann Codes and their Parameters}
\author[IC]{Ilaria Cardinali}
\ead{ilaria.cardinali@unisi.it}
\address[IC]{Department of Information Engineering and Mathematics, University of Siena,
Via Roma 56, I-53100, Siena, Italy}
\author[LG]{Luca Giuzzi\corref{cor1}}
\ead{luca.giuzzi@unibs.it}
\address[LG]{DICATAM - Section of Mathematics,
University of Brescia,
Via Branze 43, I-25123, Brescia, Italy}
\cortext[cor1]{Corresponding author.}

\begin{abstract}
  In this paper we introduce and study line Hermitian Grassmann codes as those
  subcodes of the Grassmann codes associated to the $2$-Grassmannian
  of a Hermitian polar space defined over a finite field.
  In particular, we determine the parameters and characterize the  words of minimum weight.
\end{abstract}
\end{frontmatter}

\section{Introduction}\label{sec1}
Let $V:=V(K,q)$ be a vector space of dimension $K$ over a finite field $\FF_q$
and $\Omega$ a projective system of $\PG(V)$, i.e. a set of $N$ distinct points in $\PG(V)$ such that $\dim\langle\Omega\rangle=\dim(V)$.
A projective code $\cC(\Omega)$ induced by $\Omega$ is a $[N,K]$-linear
code admitting a generator matrix $G$ whose columns are
vector representatives of the points of $\Omega$; see~\cite{tvn}.
There is a well-known relationship between the maximum number of
points of $\Omega$ lying in a hyperplane of $\PG(V)$
and the minimum Hamming distance $d_{\min}$ of $\cC(\Omega)$, namely
\[ d_{\min}=N-\max_{{\begin{subarray}{l}
                  \Pi\leq \PG(V)\\
                   \cod(\Pi)=1
                   \end{subarray}}}\left|\Pi\cap\Omega\right|. \]
Interesting cases arise when $\Omega$ is the point-set of a
Grassmann variety.
 The associated codes $\cC(\Omega)$ are called \emph{Grassman codes} and have
been extensively studied, see e.g.~\cite{R1,R2,R3,N96,GPP2009,GK2013,KP13}.

In~\cite{IL13}, we started investigating some projective codes arising
from subgeometries of the Grassmann variety associated
 to orthogonal and symplectic $k$-Grass\-mann\-ians.
 We called such codes respectively \emph{orthogonal} \cite{IL13,ILP14,IL17,IL16}
 and \emph{symplectic  Grassman codes} \cite{IL15,IL17}.
 In the cases of line  orthogonal and symplectic Grassmann codes, i.e. for $k=2$, we determined
 all the parameters; see~\cite{ILP14},~\cite{IL16} and~\cite{IL15}. For both these families
we also proposed in \cite{IL17} an efficient encoding algorithm, based
on the techniques of enumerative coding introduced in~\cite{Cover}.

In this paper we define \emph{line Hermitian Grassmann codes} as the
projective codes  defined by the projective system
of the points of the image under the
Pl\"ucker embeddings of line Hermitian Grassmannians and determine their parameters.
We refer the reader to Section~\ref{Hermitian grassmannian}
for the definition and properties of Hermitian  Grassmannians.
\medskip

\noindent {\bf Main Theorem.} {\it
 A line Hermitian  Grassmann code defined by a non--degenerate
Hermitian form on a vector space $V(m,q^2)$ is a $[N,K,d_{\min}]$-linear code where
   \[ N=\frac{(q^m+(-1)^{m-1})(q^{m-1}-(-1)^{m-1})(q^{m-2}+(-1)^{m-3})(q^{m-3}-(-1)^{m-3})}{(q^{2}-1)^2(q^{2}+1)};\]
   \[K={m\choose 2};\]
  \[d_{\min}=\begin{cases}\displaystyle
    q^{4m-12}-q^{2m-6} & \text{ if $m=4,6$ } \\
    q^{4m-12} & \text{ if $m\geq 8$ is even} \\
    q^{4m-12}-q^{3m-9} & \text{ if $m$ is odd.}
  \end{cases}.
\]}
\medskip

As a  byproduct of the  proof of the Main Theorem,
we obtain a characterization of the words of minimum weight for any $m$ and $q$, except for $(m,q)=(5,2)$, see Corollaries~\ref{codd} and \ref{ceven}.

In a forthcoming paper \cite{IL-HE} we plan to describe and discuss
algorithms for implementing encoding, decoding and error
correction for line Hermitian Grassmann codes in the same spirit
 of \cite{IL17}.

\subsection{Organization of the paper}
In Section \ref{preliminaries} we recall
some preliminaries and set our notation.
In particular, in Section~\ref{Hermitian grassmannian} some basic notions about projective codes,
Hermitian Grassmannians and their Pl\"ucker embeddings are recalled,
while
in Section~\ref{weight formula} we recall a formula for
estimating the weight of codewords for Grassmann codes.
The same formula appears also in \cite{ILP14}, but we now offer a
much simplified and shorter proof.
Section~\ref{s3} is dedicated to the proof of our main result,
by determining the minimum weight of line Hermitian
Grassmann codes and, contextually, obtaining a description of
the words of minimum weight in geometric terms.
In particular, in Section~\ref{westim} we provide bounds on the values of the weights
given by the formula of Section~\ref{weight formula} and in Sections~\ref{odd_m} and~\ref{lev0} we investigate in detail the minimum weight in the cases where the hosting space
has odd or even dimension.

\section{Preliminaries}\label{preliminaries}
\subsection{Hermitian Grassmannians and their embeddings}\label{Hermitian grassmannian}
There is an extensive literature on the properties of Hermitian
varieties over finite fields; for the basic notions as well as proofs
for the counting formulas we use, we refer  to the
monograph~\cite{S65} as well as well as to the survey \cite{BC66}; see also
\cite[Chapter 2]{HT16}.  We warn the reader that we choose to uniformly
use  vector dimension in all statements throughout this paper.

Given any  $m$-dimensional vector space $V:=V(m, \KK)$  over a field $\KK$ and $k\in \{1,\dots, m-1\}$,
let $\cG_{m,k}$  be the  $k$-Grass\-mann\-ian of the projective
space $\PG(V)$, that is the point--line
geometry whose points are the $k$-dimensional subspaces of $V$
and whose lines are the sets
\[ \ell_{W,T}:=\{ X: W\leq X\leq T, \dim X=k \} \]
with $\dim W=k-1$ and  $\dim T=k+1$.

When we want to stress on the role of the vector space $V$ rather than its dimension $m$, we shall
write $\cG_k(V)$ instead of $\cG_{m,k}.$
In general, the points of a projective space $\PG(V)$ will be denoted by $[u]$, where $u\in V$
is a non-zero vector. For any $X\subseteq V$, we shall also write $[X]:=\{ [x] \colon x\in \langle X\rangle\}$.

Let
$e_{k}:\cG_{m,k}\to\PG(\bigwedge^kV)$ be the Pl\"ucker (or Grassmann) embedding of $\cG_{m,k}$, which maps an arbitrary
$k$--dimensional subspace $X=\langle v_1,v_2,\ldots,v_k\rangle$ of
$V$ to the point $e_k(X):=[v_1\wedge v_2\wedge\cdots\wedge v_k]$ of $\PG(\bigwedge^kV)$.
Note that lines of $\cG_{m,k}$ are mapped onto (projective) lines of $\PG(\bigwedge^kV)$.
The dimension $\dim(e_k)$
of the embedding is defined to be the vector dimension of the subspace spanned by its image.
It is well known that $\dim(e_k)={m\choose k}.$

The image $e_{k}(\cG_{m,k})$ of the Pl\"ucker embedding is a projective variety
of $\PG(\bigwedge^kV)$,
called \emph{Grassmann variety} and denoted by $\GG(m,k)$.

By Chow's theorem \cite{C49}, the semilinear automorphism group stabilizing the
variety $\GG(m,k)$ is the projective general semilinear group
${\mathrm P}{\Gamma}{\mathrm L}(m,\KK)$ unless $k=m/2$, in which case it is
${\mathrm P}{\Gamma}{\mathrm L}(m,\KK)\rtimes{\mathbb Z}_2$.
This is also the permutation automorphism group of the induced code; see \cite{GK2013}.

In~\cite{ILP17} we introduced the notion of \emph{transparent embedding} $e$ of a point-line geometry $\Delta$, as a way to clarify the relationship between the automorphisms of $\Delta$ and the automorphisms of its image $\Omega:=e(\Delta)$ (and, consequently, also the automorphisms of the codes  $\cC(\Omega)$). A projective embedding $e:\Delta\to\PG(W)$ where $W=\langle \Omega\rangle$, is called \emph{(fully) transparent} when the pre-image of every line contained (as a point-set) in
$\Omega$ is actually a line of $\Delta$.
When an embedding is homogeneous and transparent, the collineations
of $\PG(W)$ stabilizing $\Omega$ lift to automorphisms of $\Delta$ and, conversely, every
automorphism of $\Delta$ corresponds to a collineation of $\PG(W)$ stabilizing $\Omega$.
So, under this assumption, it is possible to easily describe the relationship between
the groups which are involved.
In particular, the Grassmann embedding $e_k:\cG_{m,k}\to\GG(m,k)$ is transparent.

Assume henceforth $\KK=\FF_{q^2}$, so $V$ is an $m$-dimensional vector space defined over a finite field  of order $q^2.$ Suppose that $V$
is equipped with a non-degenerate Hermitian form $\eta$ of Witt index $n$ (hence $m=2n+1$ or $m=2n$).

The Hermitian $k$-Grassmannian $\cH_{n,k}$ induced by $\eta$
is defined for $k=1,\ldots,n$ as the subgeometry of $\cG_{m,k}$ having as
points the totally $\eta$--isotropic
subspaces of $V$ of dimension $k$ and as lines
\begin{itemize}
\item for $k< n$, the sets of the form
   \[ \ell_{W,T}:=\{ X: W\leq X\leq T, \dim X=k \} \]
   with $T$ totally $\eta$--isotropic and $\dim W=k-1$, $\dim T=k+1$.
 \item for $k=n$, the sets of the form
   \[ \ell_{W}:=\{ X: W\leq X, \dim X=n, X\,\, {\text{totally $\eta$--isotropic}}\} \]
   with $\dim W=n-1$, $W$ totally $\eta$--isotropic.
\end{itemize}
If $k=1$, $\cH_{n,1}$ indicates a Hermitian polar space of rank $n$ and if $k=n$, $\cH_{n,n}$ is usually called \emph{Hermitian dual polar space of rank $n$}.

Let $\varepsilon_{n,k}:=e_{k}|_{\cH_{n,k}}$ be the restriction of the Pl\"ucker embedding $e_{k}$ to the Hermitian $k$-Grassmannian $\cH_{n,k}.$ The map $\varepsilon_{n,k}$ is an embedding of $\cH_{n,k}$ called
\emph{Pl\"ucker (or Grassmann) embedding} of $\cH_{n,k}$; its dimension is proved to be $\dim(\varepsilon_{n,k})={{\dim(V)}\choose k}$ if $\dim(V)$ is even and $k$ arbitrary by Blok and Cooperstein~\cite{BC2012} and for $\dim(V)$ arbitrary and $k=2$ by Cardinali and Pasini~\cite{IP14}.

Put $\HH_{n,k}:=\varepsilon_{n,k}(\cH_{n,k})=\{\varepsilon_{n,k}(X)\colon X {\text{ point of }} \cH_{n,k}\}$. Then $\HH_{n,k}$ is a projective system of $\PG(\bigwedge^kV).$

Note that if $k=2$ and $n>2$ then $\varepsilon_{n,2}$ maps lines of $\cH_{n,2}$ onto projective lines of $\PG(\bigwedge^2 V)$, independently from the parity of $\dim(V)$, i.e. the embedding is
projective. Otherwise, if $n=k=2$ and $m=\dim(V)=5$ then the lines of $\cH_{2,2}$ are mapped onto Hermitian curves, while if $m=\dim(V)=4$ then lines of $\cH_{2,2}$ are mapped onto Baer sublines of $\PG(\bigwedge^2 V)$. In the latter case $\HH_{2,2}\cong Q^-(5,q)$ is
contained in a proper subgeometry of $\PG(\bigwedge^2V)$ defined over $\FF_q$.
We observe that for $k=2$, $\dim(V)=4$ or $\dim(V)>5$
the embeddings $\varepsilon_{n,2}$ are always transparent; see
\cite{ILP17}.

We will denote by $\cC(\HH_{n,k})$ the projective code arising from vector representatives of the elements of $\HH_{n,k}$, as explained at
the beginning of the Introduction.

The following theorem is a consequence of the transparency of
the embedding $\varepsilon_{n,2}$ and the description of the monomial automorphism group of
projective codes; see \cite{GK2013}.
\begin{theorem}
  The monomial automorphism group of the codes $\cC(\HH_{n,2})$ is
  $\mathrm{P}\Gamma\mathrm{U}(m,q)$ for $m>5$.
  For $m=4$, we have $\HH_{2,2}\cong Q^-(5,q)$; so the monomial automorphism group
  of the code is isomorphic to $GO^-(5,q)$.
\end{theorem}
Clearly, the length of $\cC(\HH_{n,k})$ is the number of points of a
Hermitian $k$-Grass\-mann\-ian $\cH_{n,k}$ and the dimension of $\cC(\HH_{n,k})$ is the dimension of the embedding $\varepsilon_k.$

From here on we shall focus on the minimum distance $d_{\min}$  of a line Hermitian code, i.e. $k=2$.
There is a geometrical way to read the minimum distance of $\cC(\HH_{n,2})$: since any
codeword of $\cC(\HH_{n,2})$ corresponds to a bilinear alternating form on $V$,
it can be easily seen that the minimum distance of $\cC(\HH_{n,2})$ is precisely the length of $\cC(\HH_{n,2})$ minus the maximum number of lines which are simultaneously totally $\eta$-isotropic
for the given Hermitian form $\eta$ defining $\cH_{n,2}$ and totally $\varphi$-isotropic for a (possibly degenerate) bilinear alternating form $\varphi$ on $V.$

\subsubsection{Notation}
Since the cases $\dim(V)$ even and $\dim(V)$ odd behave differently, it will
be sometimes useful to adopt the following notation. We will write $\cH_{n,k}^{odd}$ for a Hermitian $k$-Grassmannian in the case $\dim(V)=2n+1$ and $\cH_{n,k}^{even}$ for a Hermitian $k$-Grassmannian in the case $\dim(V)=2n.$ Accordingly,  $\mu_n^{odd}$ is the number of points of $\cH_{n,1}^{odd}$ and $\mu_n^{even}$ is the number of points of $\cH_{n,1}^{even}.$

For $k=2$, the
number of points of $\cH_{n,2}^{odd}$ is the length $N^{odd}$
of $\cC(\HH_{n,2}^{odd})$:
\begin{equation}\label{lenght odd}
N^{odd}= \frac{\mu_{n-1}^{odd}\cdot \mu_{n}^{odd}}{q^2+1}\,\,\,\,\,\,\mbox{where}\,\,\,\,\,\, \mu_n^{odd}:=\frac{(q^{2n+1}+1)(q^{2n}-1)}{(q^2-1)}.
 \end{equation}

 Analogously,
the number of points of $\cH_{n,2}^{even}$ is the length $N^{even}$ of $\cC(\HH_{n,2}^{even})$:
\begin{equation}\label{lenght even}
N^{even}= \frac{\mu_{n-1}^{even}\cdot \mu_{n}^{even}}{q^2+1}\,\,\,\,\,\, \mbox{where}\,\,\,\,\,\, \mu_n^{even}:=\frac{(q^{2n-1}+1)(q^{2n}-1)}{(q^2-1)}.
 \end{equation}
 Equations~\eqref{lenght odd} and~\eqref{lenght even} together with the results from~\cite{BC2012,IP14} on the dimension of the Grassmann embedding of a line Hermitian Grassmannian prove the first claims of
 the Main Theorem about the length and the dimension of the code.

 When we do not want to explicitly focus on the Witt index of $\eta$ but we prefer to stress on $\dim(V)=m$ regardless of its parity, we write $\cH_{m,k}$ for the Hermitian $k$-Grassmannian  defined by $\eta$ and  $\varepsilon_{m,k}$ for its Pl\"{u}cker embedding; we also put
 $\varepsilon_{m,k}(\cH_{m,k})=\HH_{m,k}$. Clearly, if $m$ is odd (i.e. $m=2n+1$) then the symbols $\cH_{m,k}$ and $\cH_{n,k}^{odd}$ have the same meaning and analogously,  if $m$ is even (i.e. $m=2n$), the symbols $\cH_{m,k}$ and $\cH_{n,k}^{even}$. Accordingly,
\[ \mu_m=\begin{cases}
    \mu^{even}_{m/2} & \text{ if $m$ is even } \\
    \mu^{odd}_{(m-1)/2}   & \text{ if $m$ is odd.}
    \end{cases} \]
  For simplicity of notation, we shall always write $\cH_m$ for the point-set of $\cH_{m,1}$.

\subsection{A recursive weight formula for Grassmann and polar Grassmann codes}
  \label{weight formula}
Denote by $V^*$ the dual of a vector space $V.$
It is well known that 
$(\bigwedge^kV)^*\cong\bigwedge^kV^*$ and that the linear functionals belonging to $(\bigwedge^kV)^*$
correspond exactly to $k$-linear alternating forms defined on $V$.
More in detail, given $\varphi\in\bigwedge^kV^*$, we have that
\[ \varphi^*(v_1,\ldots,v_k):=\varphi(v_1\wedge v_2\wedge\cdots\wedge v_k) \]
is a $k$-linear alternating form on $V$. Conversely, given any $k$-linear alternating
form $\varphi^*\colon V^k\to\FF_q$, there is a unique element $\varphi\in (\bigwedge^kV)^*$
such that
\[ \varphi(v_1\wedge\ldots\wedge v_k):=\varphi^*(v_1,\ldots,v_k) \]
for any $v_1,\ldots,v_k\in V$.
Observe that, given a point $[u]=[ v_1\wedge v_2\wedge\cdots\wedge v_k]\in\PG(\bigwedge^kV)$,
we have
$\varphi(u)=0$
if and only if all the $k$-tuples of elements of the vector space $U:=\langle v_1,\ldots,v_k\rangle$ are killed by  $\varphi^*$.
With a slight abuse of notation, in the remainder of this paper,
we shall use the same symbol $\varphi$ for both the linear functional and the
related $k$-alternating form.

Suppose $\{X_1,\ldots,X_N\}$ is a set of $k$-spaces of $V$ and consider the projective
system
$\Omega=\{[\omega_1],\ldots,[\omega_N]\}$ of $\PG(\bigwedge^kV)$ with $[\omega_i]:=e_k(X_i)$, $1\leq i\leq N$,
where $e_k$ is the Pl\"{u}cker embedding of $\cG_k(V)$.
Put $W:=\langle\Omega\rangle$
and let
\[ \cN(\Omega):=\{ \varphi\in\bigwedge^kV^*\colon \varphi|_{\Omega}\equiv 0 \} \]
be the annihilator of the set $\Omega$; clearly $\cN(\Omega)=\cN(W)$.
There exists a correspondence between the elements of
 $(\bigwedge^kV^*)/\cN(\Omega)\cong W^*$ and the codewords of $\cC(\Omega),$ where $\cC(\Omega)$ is the linear code associated to $\Omega$.
Indeed, given any $\varphi\in W^*,$
the codeword $c_{\varphi}$ corresponding to $\varphi$ is
\[ c_{\varphi}:=(\varphi(\omega_1),\ldots,\varphi(\omega_N)). \]
As $\Omega$ spans $W$ it is immediate to see that $c_{\varphi}=c_{\psi}$ if and only if
$\varphi-\psi\in\cN(\Omega)$, that is $\varphi=\psi$ as elements of $W^*$.

We define the weight $\wt(\varphi)$ of $\varphi$ as the weight of the
codeword $c_{\varphi}$
\begin{equation}\label{peso funzionale}
 \wt(\varphi):=\wt(c_{\varphi})=|\{ [\omega]\in\Omega \colon \varphi(\omega)\neq 0 \}|.
\end{equation}
When $\varphi$ is a non-null linear functional in $\bigwedge^kV^*$, its kernel determines
a hyperplane $[\Pi_{\varphi}]$ of $\PG(\bigwedge^kV)$;
hence Equation~\eqref{peso funzionale} says that the weight of a
non-zero codeword $c_{\varphi}$ is the number of points of the projective system $\Omega$ not lying on the hyperplane
$[\Pi_{\varphi}]$.

For linear codes the minimum distance $d_{\min}$ is the minimum of the weights of the non-zero codewords; so,
in order to obtain $d_{\min}$ for $\cC(\Omega)$ we need to determine the
maximum number of $k$--spaces of $V$ mapped by the Pl\"{u}cker embedding to $\Omega$
that are also $\varphi$--totally isotropic,
as  $\varphi$ is an arbitrary $k$--linear alternating form which is not identically null
on the elements of $\Omega$.

In~\cite[Lemma 2.2]{ILP14} we proved a condition relating the
weight of a codeword $\varphi$ of a $k$-(polar) Grassmann code with
the weight of codewords of suitable $(k-1)$-(polar) Grassmann codes.
Since that result will be useful also for the present paper,
we recall it in Lemma~\ref{formulona}, providing a much shorter
and easier proof.

\medskip
In order to properly state Lemma~\ref{formulona} in its more general form, so that it can be applied to obtain the weights of codes associated to arbitrary subsets of the Grassmann variety (and not only to polar Grassmann varieties), we need to set some further notation; as a consequence, the remainder of this section is unavoidably quite technical. 
In any case, we warn the reader that to obtain the weights  of a line Hermitian code one can  
just use the arguments of Remark~\ref{llre} in order to derive Equation~\eqref{ka-1} directly.

Given any vector $u\in V$ and a $k$-linear alternating
form $\varphi$, define $u\bigwedge^{k-2}V:=\{u\wedge y: y\in\bigwedge^{k-2} V\}\subseteq \bigwedge^{k-1} V$.
Put $V_u:=V/\langle u\rangle$. Clearly, for any $y\in u\bigwedge^{k-2}V$ we have $\varphi(u\wedge y)=0$.
We can now define the functional ${\bar{\varphi}}_u\in (\bigwedge^{k-1}V_u)^*\cong((\bigwedge^{k-1}V)/(u\bigwedge^{k-2}V))^*$ by
the clause
\[{\bar{\varphi}}_u:\left\{\begin{array}{l}
\bigwedge^{k-1}V_u \rightarrow \KK \\
x+(u\bigwedge^{k-2}V)\mapsto\varphi(u\wedge x)
\end{array}\right.\]
where $x\in \bigwedge^{k-1}V.$
The functional ${\bar{\varphi}}_u$ is well defined  and it can naturally be regarded as a $(k-1)$-linear
alternating form on the quotient $V_u$ of $V.$
Also observe that $\wt({\bar{\varphi}}_u)=\wt({\bar{\varphi}}_{\alpha u})$ for any non-zero scalar $\alpha$,
so the expression $\wt(\bar{\varphi}_{[u]}):=\wt(\bar{\varphi}_u)$ is well defined. Let
\[\Delta:=\{X_i:=e_k^{-1}([\omega_i])\colon [\omega_i]\in \Omega\}\]
be the set of $k$-spaces of $V$ mapped by the Pl\"ucker embedding to $\Omega$ and let
\[\Delta_u:=\{ X/\langle u\rangle: u\in X, X\in\Delta \} \text{ and }u^{\Delta}:=\langle X\in\Delta: u\in X\rangle.\]
Since $V_u^{\Delta}\leq V_u$ with $V_u^{\Delta}:=u^{\Delta}/\langle u\rangle$, we can consider
the restriction
\begin{equation}\label{phi_u}
\varphi_u:={\bar{\varphi}}_u|_{\bigwedge^{k-1}V_u^{\Delta}}
\end{equation}
of the functional
${\bar{\varphi}}_u$ to the space $\bigwedge^{k-1}V_u^{\Delta}.$

Note that when writing $\varphi_u$, we are implicitly assuming that $u$ belongs to one of the elements in $\Delta.$ We have $\wt({\bar{\varphi}}_u)=\wt(\varphi_u)$ because all points of $\Delta_u$ are,
by construction, contained in $V_u^{\Delta}$. Hence, $\Omega_u:=e_{k-1}(\Delta_u)=\{e_{k-1}(X/\langle u\rangle)\colon X/\langle u\rangle\in \Delta_u\}$ is a projective system of $\PG(\bigwedge^{k-1}V_u^{\Delta}).$
 The form $\varphi_u$
can be regarded as a codeword of the $(k-1)$-Grassmann code $\cC(\Omega_u)$ defined by the image $\Omega_u$ of $\Delta_u$ under the Pl\"{u}cker embedding $e_{k-1}$ of $\cG_{k-1}(V_u^{\Delta})$.

\begin{lemma}\label{formulona}
Let $V$ be a vector space over $\FF_q$, $\Omega=\{[\omega_i]\}_{i=1}^N$ a projective system of $\PG(\bigwedge^kV)$ and $\Delta=\{e_k^{-1}([\omega])\colon [\omega]\in \Omega\}$ where $e_k$ is the Pl\"{u}cker embedding of $V.$
Suppose $\varphi:\bigwedge^kV\to\KK$. Then
  \begin{equation}
  \label{ew0}
  \wt(\varphi)=\frac{q-1}{q^k-1}\sum_{\begin{subarray}{c}
      [u]\in\PG(V) \colon \\
      [u]\in [X], X\in\Delta
  \end{subarray}}\wt(\varphi_{[u]}).
\end{equation}
\end{lemma}
\begin{proof}
  By Equation~\eqref{peso funzionale},
  $\wt(\varphi)$ is the number of $k$-spaces of $V$ mapped to $\Omega$ and not killed by $\varphi.$  For any point $[u]$ of $\PG(V)$ such that $u$ is a vector in $X_i:=e_k^{-1}([\omega_i])$
  with $[\omega_i]\in \Omega$, the number of $k$-spaces through $[u]$ not killed by
  $\varphi$ is $\wt(\varphi_{[u]}):=\wt(\varphi_{u})$ (see Equations~\eqref{peso funzionale} and \eqref{phi_u}). Since any
  projective space $[X_i]$ with $X_i\in\Delta$
  contains $(q^{k}-1)/(q-1)$ points, the formula follows.
\end{proof}
Since each projective point  corresponds to $q-1$ non-zero vectors,
when we sum over the \emph{vectors} contained in
$X_i\in\Delta$ rather than over  projective points $[u]\in [X_i]$, Formula~\eqref{ew0}
reads as
\begin{equation}\label{ew0vettori}
  \wt(\varphi)=\frac{1}{q^k-1}\sum_{\begin{subarray}{c}
       u\in X_i\in\Delta
  \end{subarray}}\wt(\varphi_u).
\end{equation}

Even if Equations~\eqref{ew0} and~\eqref{ew0vettori} are equivalent, in this paper
we find it convenient to use more often Equation~\eqref{ew0vettori} then~\eqref{ew0}.

\begin{remark}
  \label{llre}
  For the purposes of the present paper we observe that $\wt(\varphi_u)$ could also be defined
  just as the number of $(k-1)$-subspaces of $V_u=V/\langle u\rangle$
  which are both isotropic with respect to both the
  polarities $\perp_{\eta}^u$ and $\perp_{\varphi}^{u}$
  induced by respectively $\perp_{\eta}$ and $\perp_{\varphi}$ in $V_u$ by the clauses
  \[ (x+\langle u\rangle)\perp_{\eta}^u (y+\langle u\rangle)\Leftrightarrow
    x\perp_{\eta} y; \]
   \[ (x+\langle u\rangle)\perp_{\varphi}^u (y+\langle u\rangle)\Leftrightarrow
    x\perp_{\varphi} y. \]
Later on, in Equation~\eqref{ka-1} we will rewrite Equation~\eqref{ew0vettori} for the special case of line Hermitian Grassmannians.

\end{remark}

\section{Weights for Hermitian Line Grassmann codes}
\label{s3}
In this section we shall always assume $V:=V(m,q^2)$ to be an $m$-dimensional vector space over
the finite field $\FF_{q^2}$, regardless the parity of $m$, $\eta$ to be a non-degenerate Hermitian form on $V$ with
$\cH_{m}$ the (non-degenerate) Hermitian polar space associated to $\eta$ and $\Delta:=\cH_{m,2}$ to be
the set of lines of $\cH_{m}$, i.e. the set of totally $\eta$-isotropic lines of $\PG(V).$
Since we clearly consider only the cases for which $\Delta$ is non-empty, we have $m\geq4$.

\subsection{Estimates}
\label{westim}
We start by explicitly rewriting Equation~\eqref{ew0vettori} for the case $k=2$, i.e.  for line Hermitian Grassmannian codes.
According to the notation introduced above we have $\varphi\in\bigwedge^2V^*$ and  $\Omega:=\{\varepsilon_{m,2}(\ell)\colon \ell\in \Delta\}.$
For any $\varphi\in \bigwedge^2V^*$ and
for $u\in V$, put $u^{\perp_{\eta}}=\{y\colon \eta(x,y)=0\}$ and  $u^{\perp_{\varphi}}=\{y\colon \varphi(x,y)=0\}.$
Observe that $[u]\in \cH_m$, is equivalent to $u\in u^{\perp_{\eta}}$; thus,
$u^{\perp_{\eta}}$ corresponds precisely to the set $u^{\Delta}$ defined in Section~\ref{weight formula}.
Explicitly, Equation~\eqref{phi_u} can be written as:
\begin{equation}\label{phi_u-tilde}
\varphi_u:\left\{\begin{array}{l}
u^{\perp_{\eta}}/\langle u\rangle\to\FF_{q^2} \\
{\varphi}_u(x+\langle u\rangle)=\varphi(u \wedge x) (=\varphi(u,x)).
\end{array}\right.
\end{equation}
The function $\varphi_u$ can be regarded as a linear functional on $u^{\perp_{\eta}}/\langle u\rangle$. Its kernel $\ker(\varphi_u)=(u^{\perp_{\varphi}}/\langle u\rangle) \cap (u^{\perp_{\eta}}/\langle u\rangle)$ either is the whole $u^{\perp_{\eta}}/\langle u\rangle$ or
it is a subspace $\Pi_u$ inducing a hyperplane $[\Pi_u]$ of $\PG(u^{\perp_{\eta}}/\langle u\rangle).$

Note that since $\eta(u,x)=0$ for all $x\in u^{\perp_{\eta}}$,
the vector space $u^{\perp_{\eta}}/\langle u\rangle$  is naturally  endowed with the Hermitian
form $\eta_u:(x+\langle u\rangle, y+\langle u\rangle)\to\eta(x,y)$ and  $\dim (u^{\perp_{\eta}}/\langle u\rangle)=\dim(V)-2.$  It is well known that
the set of all totally singular vectors
for $\eta_u$ defines (the point set of) a non-degenerate Hermitian polar space $\cH_{m-2}$ embedded in $\PG(u^{\perp_{\eta}}/\langle u\rangle).$

We shall now apply Equation~\eqref{ew0} to the codewords of the line Hermitian Grassmann code $\cC(\HH_{m,2}).$ To this aim, we rewrite it as
\begin{equation}
\label{ka-1-1}
\wt(\varphi)=  \frac{q-1}{q^4-1}\sum_{\begin{subarray}{c}
    [u]\in \ell\in \cH_{m,2}\\
    \end{subarray}}\!\!\!\!\!\!\wt(\varphi_u)=\,\frac{q-1}{q^4-1}\sum_{\begin{subarray}{c}
    [u]\in \cH_{m}\\
    \end{subarray}}\wt({\varphi}_u).
\end{equation}
Similarly, when considering vectors, Equation~\eqref{ew0vettori} can be rewritten as
\begin{equation}
\label{ka-1}
\wt(\varphi)=  \frac{1}{q^4-1}\sum_{\begin{subarray}{c}
    u\in V:
    [u]\in \cH_m\\
    \end{subarray}}\wt({\varphi}_u).
\end{equation}


Let $u$ be a vector such that $[u]\in\cH_{m}$.
By Equation~\eqref{peso funzionale} in Section~\ref{weight formula}, $\wt(\varphi_u)$ is the number of $\eta$-isotropic lines $\ell=[ v_1,v_2]$ of $\PG(V)$ through $[u]$ such that $\varphi(v_1,v_2)\neq0$ or, equivalently, working in the setting $u^{\perp_{\eta}}/\langle u\rangle$, $\wt(\varphi_u)$ is the number of points contained in the hyperplane $[\Pi_u]$ not lying on $\cH_{m-2}.$  The hyperplane $[\Pi_u]$ can be either secant (i.e. meeting in a non-degenerate variety) or tangent to $\cH_{m-2}$. Recall that all secant sections of $\cH_{m-2}$ are projectively equivalent to a Hermitian polar space $\cH_{m-3}$ embedded
in $\PG(m-3,q^2)$. So,
we have the following three possibilities:
\begin{enumerate}[a)]
\item\label{cA} $[\Pi_u]\cap \cH_{m-2}=\cH_{m-2}$ if $\ker(\varphi_u)\cong u^{\perp_{\eta}}/\langle u\rangle$;
\item\label{cB} $[\Pi_u]\cap \cH_{m-2}=\cH_{m-3}$ if $[\Pi_u]$ is a secant hyperplane to $\cH_{m-2}$;
\item\label{cC} $[\Pi_u]\cap \cH_{m-2}=[u]\cH_{m-4}$ if $[\Pi_u]$ is a hyperplane tangent to $\cH_{m-2}$, where $[u]\cH_{m-4}$ is a cone with
  vertex the point $[u]$ and basis a non-degenerate Hermitian polar space $\cH_{m-4}.$
\end{enumerate}
Put \begin{equation}\label{mu_m}
\mu_m:=|\cH_m|=\frac{(q^m+(-1)^{m-1})(q^{m-1}- (-1)^{m-1})}{(q^2-1)}
\end{equation}
for the number of points of $\cH_m$.
By convention, we put $\mu_0=0$.
Three possibilities can occur for the weights of $\varphi_u$, namely
\[ \wt(\varphi_u)=\begin{cases}
    0 & \text{ in case \ref{cA})} \\
    \mu_{m-2}-\mu_{m-3}={q^{2m-7}+(-1)^{m-4}q^{m-4}} & \text{ in case \ref{cB})} \\
    \mu_{m-2}-q^2{\mu_{m-4}}-1=q^{2m-7} & \text{ in case \ref{cC})}.  \\
    \end{cases} \]
For any given form $\varphi\in (\bigwedge^2 V)^*$ write,
\begin{equation}\label{parameters}
\begin{array}{lll}
     \fA_{\varphi}&:=\{ u: [u]\in\cH_{m}, \wt({\varphi}_u)=0, u\neq\mathbf{0} \} & A:=|\fA_{\varphi}| \\
     \fB_{\varphi}&:=\{ u: [u]\in\cH_{m}, \wt({\varphi}_u)={q^{2m-7}+(-1)^{m}q^{m-4}} \} & B:=|\fB_{\varphi}| \\
     \fC_{\varphi}&:=\{ u: [u]\in\cH_{m}, \wt({\varphi}_u)=q^{2m-7}\} & C:=|\fC_{\varphi}|. \\
   \end{array}
\end{equation}

Since $u$ varies among all (totally $\eta$-singular) vectors such that $[u]\in\cH_m$, we clearly have
  $A+B+C=(q^2-1)\mu_m$, and
  Equation~\eqref{ka-1} can be rewritten as
\begin{multline}
   \label{eww}
   \wt(\varphi)=\frac{q^{2m-7}(B+C)+(-1)^{m}q^{m-4}B}{q^4-1}=\\
  \frac{(q^{2m-7}+(-1)^{m}q^{m-4})(\mu_{m}(q^2-1)-A)-(-1)^{m}q^{m-4}C}{q^4-1};
\end{multline}
thus, we can express $\wt(\varphi)$ either as a function depending on $A$ and $B$ or as a function depending on $A$ and $C$ as
\begin{multline}
  \label{ewh}
  \wt(\varphi)=\frac{q^{2m-7}}{q^2+1}\mu_{m}-
\frac{q^{2m-7}}{q^4-1}A+(-1)^{m}\frac{q^{m-4}}{q^4-1}B= \\
\vspace{20pt}
\frac{(q^{2m-7}+(-1)^{m}q^{m-4})}{q^2+1}\mu_{m}-\frac{(q^{2m-7}+(-1)^{m}q^{m-4})}{q^4-1}A-(-1)^{m}\frac{q^{m-4}}{q^4-1}C.
\end{multline}

Denote by $A_{\max}$ the maximum value $A$ might assume as $\varphi$ varies among all non-trivial
bilinear alternating forms defined on $V$. Then, by the first Equation of~\eqref{ewh} with $B=0$ and by the second Equation of~\eqref{ewh} with $C=0$ we have the following lower bounds for the minimum distance
of $\cC(\HH_{m,2})$:
\begin{equation}\label{e5}
  d_{\min}\geq\begin{cases}
    \frac{q^{2m-7}}{q^2+1}\left(\mu_{m}-\frac{1}{q^2-1}A_{\max}\right)
    & \text{if $m$ is even} \\
      \\
        \frac{q^{2m-7}-q^{m-4}}{q^2+1}\left(\mu_{m}-\frac{1}{q^2-1}A_{\max}\right)
    & \text{if $m$ is odd}.
  \end{cases}
\end{equation}
We shall determine the actual values of
$d_{\min}$ and see that the bound in~\eqref{e5} is not sharp unless $m=4,6$.
More in detail, in the remainder of this paper
we shall determine the possible values of the parameter $A$ appearing in Equation~\eqref{e5}
as a function depending on the dimension $\dim(\Rad(\varphi))$ of the radical of the form $\varphi$ and show
that in all cases the minimum weight codewords occur for $A=A_{\max}$ (but, in general, $B,C\neq 0$).
We will also characterize the minimal weight codewords.

Given a (possibly degenerate) alternating bilinear form $\varphi$ on $V$,
denote by $\Rad(\varphi)$ the radical of $\varphi$, i.e. $\Rad(\varphi)=\{x\in V\colon \varphi(x,y)=0\, \forall y\in V\}$.
Define also $f_{\varphi}:\PG(m-1,q^2)\to\PG(m-1,q^2)$  as the semilinear transformation   given by
\begin{equation}\label{f_varphi}
f_{\varphi}([x]):=[x]^{\perp\varphi\perp\eta}.
\end{equation}
It is straightforward to see that $\ker(f_{\varphi})=[\Rad(\varphi)]$.
\begin{lemma}\label{lem-perp}
Let $[u]\in \cH_m.$ Then $\varphi_u=0 \Leftrightarrow u^{\perp_{\eta}}\subseteq u^{\perp_{\varphi}}.$
\end{lemma}
\begin{proof}
Take $x\in u^{\perp_{\eta}}$ and suppose $u^{\perp_{\eta}}\subseteq u^{\perp_{\varphi}}.$ Then $\varphi(u,x)=0$, so $\varphi_u(x+\langle u\rangle)=0\,\,\forall x\in u^{\perp_{\eta}}.$
Conversely, suppose $\varphi_u$ is identically null. Then $\varphi_u(x+\langle u\rangle)=\varphi(u,x)=0\,\, \forall x\in u^{\perp_{\eta}}.$ Hence $u^{\perp_{\eta}}\subseteq u^{\perp_{\varphi}}.$
\end{proof}
By Lemma~\ref{lem-perp}, $\fA_{\varphi}=\{u\colon [u]\in \cH_m, u^{\perp_{\eta}}\subseteq u^{\perp_{\varphi}}\}
=\fA_{\varphi}^{(1)}\cup \fA_{\varphi}^{(2)}$ where
\begin{equation}\label{frakA1-frakA2}
\fA_{\varphi}^{(1)}:=\{ u\colon [u]\in \cH_m, u^{\perp_{\eta}}\subset u^{\perp_{\varphi}}\}\quad\text{ and }\quad \fA_{\varphi}^{(2)}:=\{ u\colon [u]\in \cH_m, u^{\perp_{\eta}}= u^{\perp_{\varphi}}\}.
\end{equation}

The vectors $u$ such that $u^{\perp_{\eta}}\subset u^{\perp_{\varphi}}$ are precisely those vectors for which $u^{\perp_{\varphi}}=V$, hence $\fA_{\varphi}^{(1)}=\{u: [u]\in [\Rad(\varphi)]\cap \cH_m\}.$

Let us focus now on the set $\fA_{\varphi}^{(2)}.$ Note that  $u^{\perp_{\eta}}= u^{\perp_{\varphi}}$ is equivalent to $\alpha u=u^{\perp_{\varphi}\perp_{\eta}}$ for some $0\neq\alpha\in \FF_{q^2}$, or, in terms of projective points,
$[u]=[u]^{\perp_{\varphi}\perp_{\eta}}=f_{\varphi}([u])$. Hence, $\fA_{\varphi}^{(2)}=\{ u: [u]\in \Fix(f_{\varphi})\cap\cH_m\}$ where
$\Fix(f_{\varphi})=\{ [u]: f_{\varphi}([u])=[u] \}\leq \PG(t,q),\,0\leq t\leq m-1$, is
contained in a subgeometry over $\FF_q$ of $\PG(m-1,q^2)$.
\begin{lemma}
  \label{tlemma}
Let $[u]$ be a point of $\cH_m$. The following hold.
\begin{enumerate}[a)]
\item $u\in \fA_{\varphi} \Leftrightarrow f_{\varphi}([u])=[u]$ or $u\in \Rad(\varphi).$
\item $u\in \fB_{\varphi} \Leftrightarrow f_{\varphi}([u])\neq[u]$ and  $f_{\varphi}([u])$ is a non-singular point for $\eta.$
\item $u\in \fC_{\varphi} \Leftrightarrow f_{\varphi}([u])\neq[u]$ and  $f_{\varphi}([u])$ is a singular point for $\eta.$
\end{enumerate}
\end{lemma}
\begin{proof}
  By Equations~\eqref{f_varphi} and~\eqref{frakA1-frakA2} we have  $u\in \fA_{\varphi}^{(2)}$ if and only if $[u]\in \cH_m$ and  $[u]$ is a fixed point of $f_{\varphi}$. Also, $u\in \fA_{\varphi}^{(1)}$ if and only if $[u]\in [\Rad(\varphi)]\cap \cH_m$.

Suppose $u\not\in \fA_{\varphi}.$  Then $f_{\varphi}([u])\neq[u]$; hence  $[u, f_{\varphi}([u])]$ is a line.
Since $[u]\in [u]^{\perp_{\varphi}}$, we always
have $[u]^{\perp_{\varphi}\perp_{\eta}}=f_{\varphi}([u])\in [u]^{\perp_{\eta}}$.

Suppose that the point $f_{\varphi}([u])$ is non-singular with respect to $\eta$; then
$[u]^{\perp_{\varphi}}=f_{\varphi}([u])^{\perp_{\eta}}$ meets
$[u]^{\perp_{\eta}}\cap \cH_m$ in a non-degenerate polar space not containing $f_{\varphi}([u]).$ This is equivalent to saying $u\in \fB_{\varphi}.$

In case $f_{\varphi}([u])$ is singular with respect to $\eta$ we have that $[u]^{\perp_{\eta}}\cap \cH_m\cap [u]^{\perp_{\varphi}}$ is a degenerate polar space with radical of dimension $2$, i.e. with radical the line
$[u, f_{\varphi}([u])].$ This is equivalent to saying $u\in \fC_{\varphi}.$
\end{proof}

\begin{lemma}
\label{l2}
 Suppose $\varphi$ is a non-singular alternating form. If $A=(q^{m}-1)(q+1)$ then $B=0$.
\end{lemma}
\begin{proof}
By hypothesis, $m$ is necessarily even because $\varphi$ is non-singular.
Since $A=(q^{m}-1)(q+1)$ and $\varphi$ is non-singular, $A=|\fA_{\varphi}^{(2)}|$, see~\eqref{frakA1-frakA2},
and  the semilinear transformation $f_{\varphi}$ fixes a subgeometry
 $\Fix(f_{\varphi})=\Fix(f_{\varphi})\cap\cH_m
 \cong \PG(m-1,q)$ of $\PG(V)$ of maximal dimension;
so  $f_{\varphi}^2$ is the identity (as it is a linear transformation
fixing a frame), that is to say $f_{\varphi}=f_{\varphi}^{-1}$ is involutory.
Thus,
  for each point $[p]$ we have $[p]^{\perp_\varphi\perp_\eta}=f_{\varphi}([p])=
  f_{\varphi}^{-1}([p])=[p]^{\perp_\eta\perp_{\varphi}}$, i.e. the polarities
  $\perp_{\eta}$ and $\perp_{\varphi}$ \emph{commute}.
  The transformation $f_{\varphi}$ stabilizes $\cH_m$;
  indeed,
  we have  $f_{\varphi}([x])\in\cH_m$ if and only if
    \[ [x]^{\perp_{\varphi}\perp_{\eta}}=f_{\varphi}([x])\in f_{\varphi}([x])^{\perp_\eta}=[x]^{\perp_{\varphi}\perp_{\eta}\perp_{\eta}}=[x]^{\perp_{\varphi}}, \]
    whence, applying $\perp_{\varphi}$ once more, we obtain
    \[ f_{\varphi}([x])\in\cH_m \Leftrightarrow x\in x^{\perp_{\eta}} \Leftrightarrow [x]\in\cH_m. \]
    In particular, $\forall [p]\in\cH_m$, $f_{\varphi}([p])\in\cH_m$.
    So,
    by Lemma~\ref{tlemma},
    $p\in\fC_{\varphi}\cup\fA_{\varphi}$ and, in particular, $\fB_{\varphi}=\emptyset$, i.e. $B=0$.
   \end{proof}

   Fix now a basis $E$ of $V$.
   Without loss of generality, we can assume that the matrix $H$ representing the Hermitian form
   $\eta$ with respect to $E$ is the identity matrix.
   Denote by $S$ the antisymmetric matrix  representing the
   (possibly degenerate) alternating form $\varphi$ with respect to $E$;
   recall that $\Rad(\varphi)$ is precisely the kernel  $\ker(S)$ of the matrix $S.$

   Under these assumptions, the transformation $f_{\varphi}$ can be represented as $f_{\varphi}([x]):=[S^qx^q], \, \forall x\in V.$
   Since the fixed points of a semilinear transformation of $\PG(m-1,q^2)$ are
   contained in a subgeometry  $[\Sigma_{\varphi}]\cong \PG(t,q)$ with $0\leq t\leq m-1$,
\begin{equation}\label{A-exp}
\fA_{\varphi}\subseteq\{u: [u]\in ([\Rad(\varphi)]\cap \cH_m) \cup ([{\Sigma}_{\varphi}]\cap \cH_m)\}.
\end{equation}
Put $\widetilde{\Sigma}_{\varphi}:=\FF_{q^2}\otimes {\Sigma}_{\varphi};$ $\dim(\Sigma_{\varphi})=\dim(\widetilde{\Sigma}_{\varphi})$ where $\Sigma_{\varphi}$, respectively $\widetilde{\Sigma}_{\varphi}$, is regarded as a vector space over $\FF_q$, respectively over $\FF_{q^2}.$
It is easy to see that $\Rad(\varphi)$ and $\widetilde{\Sigma}_{\varphi}$ are subspaces of $V(m,q^2)$ intersecting trivially. Since $\Rad(\varphi)=\ker(S)$, we have $\dim(\ker(S))+\dim(\widetilde{\Sigma}_{\varphi})=\dim(\ker(S))+\dim({\Sigma}_{\varphi})\leq m.$ Clearly, $\rank(S)=m-\dim(\ker(S))$, so $\dim(\widetilde{\Sigma}_{\varphi})=\dim({\Sigma}_{\varphi})\leq \rank(S)$, where $\rank(S)$ is the rank of the matrix $S.$

Put $2i:=\rank(S)$. Hence $\dim(\Rad(\varphi))=m-2i$ and $0< 2i\leq m$.
 Define
 \begin{equation}
 \label{Ai}
 A_i:=\max\{ |\fA_{\varphi}| \colon \dim(\Rad(\varphi))=m-2i \}.
 \end{equation}
 Note that if $i=0$, $\varphi$ is identically null and this gives the $\mathbf{0}$ codeword.
Clearly, by~\eqref{A-exp},
\begin{equation}
  \label{eAA}
\begin{array}{ll}
  A_i\leq &(|\Sigma_{\varphi}| +|([\Rad(\varphi)]\cap\cH_{m})|)(q^2-1)= \\
  &(q^2-1)\left(\frac{(q^{2i}-1)}{(q-1)}+|([\Rad(\varphi)]\cap\cH_{m})|\right) = \\ &(q^{2i}-1)(q+1)+|([\Rad(\varphi)]\cap\cH_{m})|(q^2-1).
  \end{array}
\end{equation}
We shall need the following elementary technical lemma.
\begin{lemma}
\label{l-sing}
Let $H$ be a non-singular matrix of order $m$ and let $t\leq m$.
  Then for
  an $(m-t)\times (m-t)$ submatrix $M$ of $H$, we have
  $m-2t\leq \rank(M)\leq (m-t)$.
\end{lemma}
\begin{proof}
  The submatrix $M$ is obtained from $H$ by deleting $t$ rows and $t$ columns.
  First delete $t$ rows. Then the rank of the $(m-t)\times m$ matrix $M'$ so obtained is
  $m-t$. If we now delete $t$ columns from $M'$ as to obtain $M$, the rank
  of $M'$ decreases by at most $t$. So, $\rank(M')-t\leq \rank(M)\leq \rank(M').$
\end{proof}

We want to explicitly determine the cardinality of $[\Rad(\varphi)]\cap\cH_m$.
The $(m-2i)$-dimensional space $[\Rad(\varphi)]$ intersects $\cH_m$ in
a (possibly) degenerate Hermitian variety.
Write $[\Rad(\varphi)]\cap\cH_m=[\Pi_t]\cH_{m-2i-t}$ where $[\Pi_t]\cH_{m-2i-t}$ is
a degenerate Hermitian variety contained in $[\Rad(\varphi)]$ with radical $[\Pi_t]$ of dimension
$t$.

By Lemma~\ref{l-sing}, $0\leq t\leq 2i$. Moreover, since $\Pi_t\subseteq \Rad(\varphi)$, $0\leq t\leq m-2i.$ (Recall also that $2i\leq m$.)
If $t=m-2i$ then $\Rad(\varphi)=\Pi_{m-2i}$ and in this case we put $[\Pi_{m-2i}]\cH_0:=[\Pi_{m-2i}].$

The following function provides the
number of points of $[\Pi_t]\cH_{m-2i-t}$ in dependence of $t$.
\begin{equation}\label{numero punti}
\begin{array}{l}
\mu_{m-2i}\colon \{0,\ldots, \min\{2i,m-2i\}\}\rightarrow \mathbb{N}\\
\mu_{m-2i}{(t)}=\begin{array}[t]{l}\displaystyle q^{2t}\mu_{m-2i-t}+\frac{q^{2t}-1}{q^2-1}=\\
       \\
                  \displaystyle{\frac{q^{2t}(q^{m-2i-t}+ (-1)^{m-2i-t-1})(q^{m-2i-t-1}- (-1)^{m-2i-t-1})+ q^{2t}-1}{q^2-1}},
	 \end{array}
\end{array}
\end{equation}
where, by convention, $\mu_{0}(0)=0$ and
$\mu_{m-2i-t}$ (for $t<m-2i$)
is the number of points of a non-degenerate Hermitian polar space $\cH_{m-2i-t}$ (see Equation~\eqref{mu_m}).

Using Equation~\eqref{numero punti}, we can rewrite Equation~\eqref{eAA} as
\begin{equation}\label{eq-A}
A_i\leq (q^{2i}-1)(q+1) + (q^2-1)\mu_{m-2i}{(t)}.
\end{equation}
\begin{lemma}\label{ll0}
For any $t\in \{0,\ldots, \min\{2i,m-2i\}\}$, we have
  $\mu_{m-2i}(t)\leq\mu_{m-2i}^{\max}$, where
  \[ \mu_{m-2i}^{\max}:=\begin{cases}
      \mu_{m-2i}(m-2i)=\frac{q^{2m-4i}-1}{q^2-1} & \text{if $i\geq m/4$}; \\
      \mu_{m-2i}{(2i)} & \text{if $i<m/4$ and $m$ is even}; \\
      \mu_{m-2i}{(2i-1)} & \text{if $i<m/4$ and $m$ is odd}.
    \end{cases} \]
\end{lemma}
\begin{proof}
  Let $H'$ be the
  matrix representing the restriction $\eta':=\eta|_{\Rad(\varphi)}$ of the Hermitian form $\eta$ to $\Rad(\varphi)$.
  Then, by Lemma~\ref{l-sing},  $m-4i\leq \rank H'\leq m-2i$.

  When $i\geq m/4$, i.e. $m-2i\leq 2i$ (and hence $0\leq t\leq m-2i$), the maximum number of points for $[\Pi_t]\cH_{m-2i-t}$ is attained for $t=m-2i$, i.e. $[\Pi_{m-2i}]\cH_{0}= [\Rad(\varphi)]=[\Pi_{m-2i}]$.
Indeed, if $m\leq 4i$ and $t=m-2i$, it is always possible to construct an antisymmetric form $\varphi$ such that $\eta'$ is the null form. This implies that  $\mu_{m-2i}(t)<\mu_{m-2i}(m-2i)=\frac{q^{2m-4i}-1}{q^2-1}$ for any $t\in \{0,\dots, \min\{m-2i\}\}$, since $[\Pi_t]\cH_{m-2i-t}\subseteq[\Rad(\varphi)].$ Hence, in this case, $\mu_{m-2i}^{\max}:=\mu_{m-2i}(m-2i)$.

  Suppose now $i<m/4$. Then, by a direct computation
  \begin{equation}
    \label{eee1}
    \mu_{m-2i}(t)-\mu_{m-2i}(t+1)=(-1)^{m-2i-t-2}q^{m-2i+t-1}.
  \end{equation}
  So,
  \begin{equation}\label{m-2i} \mu_{m-2i}(t+2)-\mu_{m-2i}(t)=(-1)^{m-t}q^{m-2i+t-1}(q-1).
  \end{equation}
  Assume $m$ even; then, by~\eqref{m-2i}, if $t$ is even,
  $\mu_{m-2i}(t+2)>\mu_{m-2i}(t)$, i.e. $\mu_{m-2i}(t)$ is a monotone increasing function
  in $t$ even.
  If $t$ is odd, then by~\eqref{m-2i},
  $\mu_{m-2i}(t+2)<\mu_{m-2i}(t)$, i.e. $\mu_{m-2i}(t)$ is a monotone decreasing function
  in $t$ odd. By~\eqref{eee1}, $\mu_{m-2i}(0)>\mu_{m-2i}(1)$.
  Recall that $0\leq t\leq 2i$; so we have
  \begin{multline*}
    \mu_{m-2i}(2i-1)<\mu_{m-2i}(2i-3)<\cdots
    <\mu_{m-2i}(1)<\\ \mu_{m-2i}(0)<
    \mu_{m-2i}(2)<\mu_{m-2i}(4)<\cdots<\mu_{m-2i}(2i).
  \end{multline*}
  In particular, the maximum value of $\mu_{m-2i}(t)$ for $i<m/4$ and $m$ even is assumed for $t=2i$, i.e. $\mu_{m-2i}^{\max}=\mu_{m-2i}(2i)$.\par

  Assume $m$ odd; then, by~\eqref{m-2i}, if $t$ is even,
  $\mu_{m-2i}(t+2)<\mu_{m-2i}(t)$, i.e. $\mu_{m-2i}(t)$ is a monotone decreasing function
  in $t$ even.
  If $t$ is odd, then by~\eqref{m-2i},
  $\mu_{m-2i}(t+2)>\mu_{m-2i}(t)$, i.e. $\mu_{m-2i}(t)$ is a monotone increasing function
  in $t$ odd. By~\eqref{eee1}, $\mu_{m-2i}(0)<\mu_{m-2i}(1)$.
  Recall that $0\leq t\leq 2i$; so we have
  \begin{multline*}
    \mu_{m-2i}(2i)<\mu_{m-2i}(2i-2)<\cdots
    <\mu_{m-2i}(2)<\\ \mu_{m-2i}(0)<
    \mu_{m-2i}(1)<\mu_{m-2i}(3)<\cdots<\mu_{m-2i}(2i-1).
  \end{multline*}
  In particular, the maximum value of $\mu_{m-2i}(t)$ for $i<m/4$ and $m$ odd is assumed for $t=2i-1$, i.e. $\mu_{m-2i}^{\max}=\mu_{m-2i}(2i-1)$.
\end{proof}
Define the function $\xi_m\colon \{1,\dots, \lfloor m/2\rfloor\}\rightarrow \mathbb{N}$
\begin{equation}
  \label{xim}
  \xi_m(i):=(q^{2i}-1)(q+1)+(q^2-1)\mu_{m-2i}^{\max},
\end{equation}
where $\mu_{m-2i}^{\max}$, introduced in Lemma~\ref{ll0}, is regarded as a function in $i$.
\begin{corollary}\label{eqr}
The following hold.
\begin{enumerate}[a)]
  \item\label{AiA} $A_i\leq\xi_m(i)$;
  \item\label{AiB}
  \begin{equation}\label{AiBd}
  d_{i}\geq\begin{cases}
    \frac{q^{2m-7}}{q^2+1}\left(\mu_{m}-\frac{1}{q^2-1}\xi_m(i)\right)
    & \text{if $m$ is even} \\
      \\
        \frac{q^{2m-7}-q^{m-4}}{q^2+1}\left(\mu_{m}-\frac{1}{q^2-1}\xi_m(i)\right)
    & \text{if $m$ is odd},
  \end{cases}
\end{equation}
    where $d_i$ is the minimum weight of the words corresponding to bilinear
    alternating
    forms $\varphi$ with $\dim(\Rad(\varphi))=m-2i$.
  \end{enumerate}
\end{corollary}
\begin{proof}
 Case~\ref{AiA}) follows from Lemma~\ref{ll0} and Equation~\eqref{eq-A}.
 Case~\ref{AiB}) follows from Equation~\eqref{ewh}, the definition~\eqref{Ai} of $A_i$ and case a).
 \end{proof}
Note that the function $\xi_m(i)$ (see~\eqref{xim}) is not monotone in $i$.
 The following lemma provides
its largest and second largest values.
  \begin{lemma}
    \label{lA2}
    If $m\neq 4,6$ the maximum value assumed by the function $\xi_m(i)$
    is attained for $i=1$.
    If $m=4,6$, then the maximum of $\xi_m(i)$ is attained for $i=m/2$.
    If $m=5$ then $\xi_5(1)=\xi_5(2).$

    The second largest value of $\xi_m(i)$ is attained for
    \[ \left\{\begin{array}{ll}
         i=1  &  \text{if $m=6$};\\
         i=3  & \text{if $m=7$}; \\
         i=4   & \text{if $m=8,9$};\\
         i=5 & \text{if $m=10$}; \\
         i=2 & \text{if $m>10$}. \\
         \end{array} \right.\]
  \end{lemma}
  \begin{proof}
    Recall that the polar line Grassmannian $\cH_{m,2}$ is non-empty only for $m\geq 4$.
\begin{itemize}
 \item  If $m=4$ then the possible values that the function $\xi_4(i)$ can assume are $\xi_4(1)$ and $\xi_4(2)$. Precisely,
    \[ \xi_4(1)=q^4+q^3+q^2-q-2<\xi_4(2)=q^5+q^4-q-1. \]
 \item If $m=5$ then $\xi_5(1)=q^5+q^4+q^2-q-2=\xi_5(2)$.
 \item If $m=6$ then the possible values that the function $\xi_6(i)$ can assume are the following
    \begin{eqnarray*}
        \xi_6(1)&=&q^7+q^6-q^5+q^3+q^2-q-2;\\
        \xi_6(2)&=&q^5+2q^4-q-2;\\
        \xi_6(3)&=&q^7+q^6-q-1.
      \end{eqnarray*}
      Hence  $\xi_6(3)>\xi_6(1)>\xi_6(2)$.
      So, for both $m=4$ and $m=6$, $\xi_m(1)<\xi_m(m/2)$ and the
      maximum value of $\xi_m(i)$ is attained for $i=m/2$.
\end{itemize}
We assume henceforth $m>6$.
For any two functions $f(x)$ and $g(x)$ we shall write
$f(x)=O_q(g(x))$ if
\[ \frac{1}{q}g(x)< f(x)< q\cdot g(x). \]
If $m$ is even and $i=1$ then
    \begin{multline}\label{m even}
      \xi_m(1)=(q^2-1)(q+1)+(q^2-1)\mu_{m-2}^{\max}=
      (q^2-1)(q+1+\mu_{m-2}{(2)})=\\
      q^{2m-5}+q^{m}-q^{m-1} +q^3+q^2-q-2=O_q(q^{2m-5}).
        \end{multline}
If $m$ is odd and $i=1$  then
    \begin{multline}\label{m odd}
      \xi_m(1)=(q^2-1)(q+1)+(q^2-1)\mu_{m-2}^{\max}=
      (q^2-1)(q+1+\mu_{m-2}(1))=\\
      q^{2m-5}+q^{m-1}-q^{m-2}+q^3+q^2-q-2=O_q(q^{2m-5}).
    \end{multline}
      If $i\geq m/4$ then $\xi_m(i):=(q^{2i}-1)(q+1)+(q^{2m-4i}-1)$. Since we always have $2i\leq m$ (so $2i\leq m\leq 4i\leq 2m,$ and clearly $i>1$), we have
      \begin{multline} \label{i>1}
            \xi_m(i)=(q^{2i}-1)(q+1)+(q^{2m-4i}-1)\leq\\
      \leq (q^{m}-1)(q+1)+(q^{m}-1)=O_q(q^{m+1}).
    \end{multline}

Hence, for $m>6$, by Equations~\eqref{m even},~\eqref{m odd},~\eqref{i>1}, we have $\xi_m(1)>\xi_m(i)$ for any  $i\geq m/4 (\geq 2)$.

We prove that also for any $1<i<\lfloor m/4\rfloor$ we have
$\xi_m(1)>\xi_m(i)$. Note that now  $m$ can be
    either even or odd.
    By Equations~\eqref{mu_m} and \eqref{numero punti}, we have that
    \[ \mu_x(t)=O_q(q^{2t}\mu_{x-t}+q^{2t-2})=O_q(q^{2t}q^{2x-3-2t}+q^{2t-2}); \]
    so $\mu_x(t)=O_q(q^{2x-3})$ for all $x\leq m$ and $t<x\leq m$, while $\mu_x(x)=O_q(q^{2x-2})$.
Hence, for $x=m-2i$,
 \[O_q(q^{2m-4i-3})<\mu_{m-2i}^{\max} < O_q(q^{2m-4i-2}),\]
  since $\mu_{m-2i}(t)>O_q(q^{2m-4i-3})$ and $\mu_{m-2i}(t)<O_q(q^{2m-4i-2})\,\,\forall t\leq m-2i.$
By Equation~\eqref{xim}, we obtain
    \[ \xi_m(i)= O_q(q^{2i+1}+q^2\mu_{m-2i}^{\max})\leq
      O_q(q^{2i+1}+q^{2m-4i}). \]
    However, for $2\leq i\leq \lfloor m/2\rfloor-1$ (henceforth also for $1<i<\lfloor m/4\rfloor$),
    \[ q^{2i+1}+q^{2m-4i}<q^m+q^{2m-8}, \]
    so
    \[ \xi_m(i)<O_q(q^{m}+q^{2m-8}). \]
    This latter value is smaller than $\xi_m(1)=O_q(q^{2m-5})$ (see Equations~\eqref{m even} and~\eqref{m odd}).
    It follows that the maximum of $\xi_m(i)$ is attained for $i=1$ for all cases $m>6$.
\par\smallskip
Assume now $i>2$. Since $O_q(q^{4m-11})<\mu_{m-4}^{\max}<O_q(q^{2m-10})$ and $\mu_{m-2i}^{\max}<O_q(q^{2m-4i-2})$ we have
$\mu_{m-4}^{\max}-\mu_{m-2i}^{\max}> O_q(q^{4m-11}-q^{2m-4i-2}).$
    Hence
      \begin{multline*} \xi_m(2)-\xi_m(i)=(q^{4}-q^{2i})(q+1)+(q^2-1)(\mu_{m-4}^{\max}-\mu_{m-2i}^{\max})\geq \\
     \geq  O_q(q^5+q^4-q^{2i+1}-q^{2i}+q^{2m-9}-q^{2m-4i}).
    \end{multline*}
    For $m/2\geq i\geq 3$, we have
    \[ q^5+q^4-q^{2i+1}-q^{2i}+q^{2m-9}-q^{2m-4i}>q^5+q^4-q^{m+1}-q^{m}+q^{2m-9}-q^{2m-12}>0, \]
    so, for $m>10$,
    \[ \xi_m(2)-\xi_m(i)>O_q(q^{2m-9})>0. \]
    A direct computation gives the following:
    \[\xi_7(1)>\xi_7(3)>\xi_7(2);\,\,\,\xi_8(1)>\xi_8(4)>\xi_8(2)>\xi_8(3); \]
   \[\xi_9(1)>\xi_9(4)>\xi_9(2)>\xi_9(3);\,\,\,\xi_{10}(1)>\xi_{10}(5)>\xi_{10}(4)>\xi_{10}(3).\]
   This completes the proof.
\end{proof}
By Corollary~\ref{eqr}  and Lemma~\ref{lA2} we have:
\begin{corollary}
  Let $\varphi$ be a form with $\dim(\Rad(\varphi))=m-2i$. Then,
  \begin{itemize}
  \item
    for $m>6$, we have $|\fA_{\varphi}|\leq A_i\leq \xi_m(1)$;
  \item
    for $m=4,6$ we have $|\fA_{\varphi}|\leq A_i\leq \xi_m(m/2)$.
  \end{itemize}
\end{corollary}
\subsection{Minimum distance of $\cH_{m,2}$ with $m$ odd}
\label{odd_m}
In this section we assume $m$ to be odd. Then the Witt index of the Hermitian form $\eta$ is  $n=(m-1)/2.$
Let $\varphi$  be an alternating form on $V.$
Recall from Equation~\eqref{ewh} that
\[
\wt(\varphi)=
\frac{(q^{2m-7}-q^{m-4})(\mu_{m}(q^2-1)-A)}{q^4-1}+\frac{q^{m-4}}{q^4-1}C.
\]
\begin{prop}\label{min wgt m odd}
There exists a bilinear alternating form $\varphi$ with $\dim(\Rad(\varphi))=m-2$ such that $\wt(\varphi)=q^{4m-12}-q^{3m-9}$ and $\wt(\varphi')\geq q^{4m-12}-q^{3m-9}$ for any other form $\varphi'$ with $\dim(\Rad(\varphi'))=m-2.$
\end{prop}
\begin{proof}
In order to determine the weight of the word of $\cH_{m,2}$ induced by the form $\varphi$, we need to determine
the number of lines of $\cH_{m,2}$ which are not totally isotropic for $\varphi$.

Take $\varphi$ with $\dim(\Rad(\varphi))=m-2.$ Then  a line $\ell$ is totally isotropic for $\varphi$ if and only if $\ell\cap[\Rad(\varphi)]\neq\emptyset.$
If $S$ denotes the matrix representing $\varphi$, we have $\rank(S)=2.$ According to the notation of
Section~\ref{westim}, $\rank(S)=2$ is equivalent to $i=1$ and $[\Rad(\varphi)]\cap \cH_m=[\Pi_t]\cH_{m-2-t}$ is a degenerate Hermitian variety with radical $[\Pi_t]$ of dimension $t.$  By Equation~\eqref{eq-A} and Lemma~\ref{ll0}, since $m>4$ is odd, the maximum number of points $\mu_{m-2}^{\max}$ of $[\Rad(\varphi)]\cap \cH_m$ is attained for $t=2i-1=1$, i.e. $\mu_{m-2}^{\max}=\mu_{m-2}(1)=q^2\mu_{m-3}+1$ (last equality comes from Equation~\eqref{numero punti}).
Hence the number of points of $\cH_m\setminus[\Rad(\varphi)]$ (see Equation~\eqref{numero punti}) is at least
\[\mu_m-\mu_{m-2}^{\max}= \mu_{m}-q^2\mu_{m-3}-1=q^{m-2}(q^{m-1}+q^{m-3}-1)=q^{2m-3}+q^{2m-5}-q^{m-2}. \]
Assume $[\Rad(\varphi)]\cap\cH_{m}=[\Pi_1]\cH_{m-3}$ and
consider a point $[p]\in\cH_{m}\setminus [\Rad(\varphi)]$.
\begin{itemize}
\item
  \framebox{Case $[p]\in [\Pi_1]^{\perp_{\eta}}$}. Then $[p]^{\perp_\eta}\cap[\Rad(\varphi)]$ is a degenerate Hermitian polar space $[\Pi_1]\cH_{m-4}$ with radical $[\Pi_1]$  of dimension $1$; so there are $(\mu_{m-2}-q^2\mu_{m-4}-1)=q^{2m-7}$ lines  through $[p]$ disjoint from $[\Rad(\varphi)]$.
The number of points $[p]$ collinear with the point $[\Pi_1]$ in $\cH_m$ but not contained
in $\cH_m\cap[\Rad(\varphi)]=[\Pi_1]\cH_{m-3}$ is $q^2(\mu_{m-2}-\mu_{m-3})=q^2(q^{2m-7}-q^{m-4}).$
\item
  \framebox{Case $[p]\not\in[\Pi_1]^{\perp_{\eta}}$}. Then $[p]^{\perp_\eta}\cap[\Rad(\varphi)]=\cH_{m-3}$, so, there are $(\mu_{m-2}-\mu_{m-3})=(q^{2m-7}-q^{m-4})$ lines
through $[p]$ which are not totally isotropic. The number of points not collinear with $[\Pi_1]$ in $\cH_m$ and not in $\cH_m\cap[\Rad(\varphi)]$ is $(\mu_{m}-q^2\mu_{m-3}-1)-q^2(\mu_{m-2}-\mu_{m-3})=(\mu_{m}-q^2\mu_{m-2}-1)=q^{2m-3}.$
\end{itemize}

So, we have that the total number of lines disjoint from $[\Rad(\varphi)]$ is
  \[ \frac{1}{q^2+1}\left(q^2(q^{m-7}-q^{m-4})\cdot q^{2m-7}+(q^{2m-7}-q^{m-4})q^{2m-3}\right), \]
  i.e.
  \[ \wt(\varphi)=\frac{1}{q^2+1}\left(q^{4m-12}-q^{3m-9}+q^{4m-10}-q^{3m-7}\right)=q^{4m-12}-q^{3m-9}. \]

  In case $[\Rad(\varphi)]\cap\cH_m\neq[\Pi_1]\cH_{m-3}$, the number of totally $\eta$-isotropic lines incident
  (either in a point or  contained in)  $[\Rad(\varphi)]$ is smaller; thus the weight of the
  word induced by $\varphi$ is larger than the value obtained above.
\end{proof}
We claim that $q^{4m-12}-q^{3m-9}$ is actually the minimum distance for $m$ odd.

As before, let $i=(\rank(S))/2$.
When $i=1$, then by Proposition~\ref{min wgt m odd}, the minimum weight of the codewords induced by $S$  is $d_1=q^{4m-12}-q^{3m-9}$.
Suppose now $i>1$;  we need to distinguish several cases according to the value of $m$.
\begin{itemize}
\item
 \framebox{$m\geq 11$}. Then by Corollary~\ref{eqr} and Lemma~\ref{lA2}
 $A_{i}\leq\xi_m(i)\leq\xi_m(2)\leq \xi_m(1).$
  By Case~\ref{AiB}) of Corollary~\ref{eqr},
  \[ d_i\geq \frac{q^{2m-7}-q^{m-4}}{q^4-1}((q^2-1)\mu_m-\xi_m(i))\geq
    \frac{q^{2m-7}-q^{m-4}}{q^4-1}((q^2-1)\mu_m-\xi_m(2)). \]
We will show that
  \[\frac{q^{2m-7}-q^{m-4}}{q^4-1}((q^2-1)\mu_m-\xi_m(2)) >q^{4m-12}-q^{3m-9}. \]
Actually, by straightforward computations, this becomes
  \[ q^{m-4}(q^{m-3}-1)(q^{2m-9}-q-1-\frac{q^{m-2}}{q^2+1})>0 \]
  which is true for all values of $q$; so $d_{\min}=d_1=q^{4m-12}-q^{3m-9}$.
\item \framebox{$m=9$}. By Lemma~\ref{lA2} we have that the second largest value of $\xi_m(i)$ is for $i=4$
  and
  $\xi_{9}(4)=q^{9}+q^8+q^2-q-2.$
This corresponds to the following bound on the minimum weight $d_i$ of codewords associated with
matrices $S$ with $2i=\rank(S)>2$ (see Case~\ref{AiB}) of Corollary~\ref{eqr}):
\[d_i>\frac{q^{11}-q^5}{q^4-1}(q^{17}-2q^9-q^2+q+1) >q^{24}-q^{18}.\]
 So, the minimum distance is attained by codewords corresponding to matrices $S$ of rank $2$
and, by Proposition~\ref{min wgt m odd},  $d_{\min}=d_1=q^{4m-12}-q^{3m-9}$.

\item \framebox{$m=7$}. By Lemma~\ref{lA2} we have that the second largest value of $\xi_m(i)$ is for $i=3$ and
  $\xi_7(3)=q^7+q^6+q^2-q-2.$
 This corresponds to the following bound on the minimum weight of codewords associated with
matrices $S$ with $2i=\rank(S)>2$:
 \[d_i>q^{16}-2q^{10}-q^5+q^4+q^3>q^{16}-q^{12}.\]
 So the minimum distance is attained by codewords corresponding to matrices $S$ of rank $2$
and, by Proposition~\ref{min wgt m odd},  $d_{\min}=d_1=q^{4m-12}-q^{3m-9}$.

\item \framebox{$m=5$}. By Lemma~\ref{lA2}, we have $\xi_5(1)=\xi_5(2)=q^5+q^4+q^2-q-2$.
We will prove that the minimum distance of $\cH_{m,2}$ is $q^8-q^6$.

 Let $\varphi$ be a (non-null) alternating bilinear form of $V(5,q^2)$ represented by a matrix $S$. The radical of $\varphi$ can have dimension $1$ or $3$, hence $\rank(S)$ is either $2$ or $4$, i.e. $i=1$ or $2.$
 By Proposition~\ref{min wgt m odd}, there exists an alternating bilinear form $\varphi$ with $\dim(\Rad(\varphi))=3$ such that $\wt(\varphi)=q^8-q^6$ and any other form $\varphi$ with $\dim(\Rad(\varphi))=3$ has weight greater then $q^8-q^6.$ So, we need to show that there
 are no alternating bilinear forms with $\dim(\Rad(\varphi))=1$ inducing words of weight
 less than $q^8-q^6$.
 Assume henceforth that $\varphi$ is an alternating bilinear form with $\dim(\Rad(\varphi))=1$ (hence $i=2$). We shall determine a lower bound $d_2$ for the weights $\wt(\varphi)$ and prove
 $\wt(\varphi)\geq d_2\geq q^8-q^6.$

By Lemma~\ref{tlemma}, $p\in\fC_{\varphi}$ if and only if $f_{\varphi}([p])\neq [p]$ and
$f_{\varphi}([p])\in\cH_5$.
If $[p]=[\Rad(\varphi)]$, then $p\in\fA_{\varphi}$.
So, suppose $[p]\neq [\Rad(\varphi)]$.
Since $[\Rad(\varphi)]=\ker(f_{\varphi})$, if $[p]\in\cH_5$ and $f_{\varphi}([p])=[x]\in\cH_5$, then $f_{\varphi}([p+\alpha \Rad({\varphi})])=[x]\in\cH_5$ for
any $\alpha\in\FF_{q^2}$.
On the other hand, for any given point $[p]\in\cH_5$, the line $[p,\Rad(\varphi)]$ meets $\cH_5$ in either $(q+1)$ or
$(q^2+1)$ points; this yields that any point in $\cH_5$ belonging to the image $f_{\varphi}(\cH_5)$ of $f_{\varphi}$ restricted to $\cH_5$ admits at least $q-1$ preimages in $\cH_5$ distinct from
itself.
Since $\dim(\Rad(\varphi))=1$, the set $\Fix(f_{\varphi})$ is contained in a $\PG(V')$ with
$\dim V'=4$. We need a preliminary lemma.
\begin{lemma}
  \label{v4c1}
  Let $f_{\varphi}:\PG(V')\to\PG(V')$ be a semilinear collineation with $\dim V'=4$.
  Then, either $\Fix(f_{\varphi})\cong\PG(3,q)$ or $|\Fix(f_{\varphi})|\leq q^2+q+2$.
\end{lemma}
\begin{proof}
  In general, $\Fix(f_{\varphi})$ is contained in a subgeometry
  $\PG(V'')$ with $V''$ a vector space over $\FF_q$ and $\dim V''=4$.
  
  Assume $\Fix(f_{\varphi})\neq\PG(V'')\cong\PG(3,q)$.
  Suppose first  that there is a vector space $W$ over $\FF_q$ with
  $\dim W=3$ such that $\PG(W)\subseteq\Fix(f_{\varphi})$.
  If $\Fix(f_{\varphi})=\PG(W)$, then we are done. Otherwise, let
  $[p]\in\Fix(f_{\varphi})\setminus\PG(W)$ and define $W':=W+\langle p\rangle_q$.
  Since $\dim_q W'=4$,
  we have that $\Fix(f_{\varphi})$ is contained in the subgeometry $\PG(W')$.
  If there is $[r]\in\Fix(f_{\varphi})\setminus\PG(W)$
  with $[r]\neq [p]$, then the subline $\ell_q:=[p,r]_q$ spanned by $[p]$ and $[r]$ meets
  the subplane $\PG(W)$ in a point $[s]\neq [r]$ which is also fixed. So $\ell_q$ is fixed
  pointwise. Take $[t]\in\PG(W)$ with $[t]\neq [s]$. The subline $[t,s]_q$ is also
  fixed pointwise, so the subplane $[p,s,t]_q\neq\PG(W)$ is also fixed pointwise.
  As $f_{\varphi}$ fixes two (hyper)planes pointwise in $\PG(W')\cong\PG(3,q)$ we
  have that $f_{\varphi}$ fixes $\PG(W')$ pointwise --- a contradiction.
  Thus, in this case $|\Fix(f_{\varphi})|\leq q^2+q+2$.

  Suppose now that $\Fix(f_{\varphi})$ does not contain a subplane isomorphic
  to $\PG(2,q)$ and that there is a subline $\ell\subseteq\Fix(f_{\varphi})$.
  If there were a subplane $\pi_q$ through $\ell$ such that
  $\Fix(f_{\varphi})\cap\pi_q$
  contains two points not on $\ell$, then (by the same argument we used in the
  case above), this subplane would have to be fixed pointwise by $f_{\varphi}$.
  This is a contradiction; so if $\Fix(f_{\varphi})$ contains a subline $\ell$,
  then there is at most one fixed point $[x]\not\in\ell$
  on any subplane through $\ell$ contained in $\PG(V'')$. So,
  $|\Fix(f_{\varphi})|\leq 2q+2$.
  
  Finally, if $\Fix(f_{\varphi})$ does not contain sublines, then $\Fix(f_{\varphi})$
  cannot contain either frames or more than $3$ points on a plane or more than $2$
  points on a line. It  follows that $|\Fix(f_{\varphi})|\leq 4$. This completes
  the proof.
\end{proof}
In light of Lemma~\ref{v4c1}
we now distinguish two subcases:
\begin{enumerate}[a)]
\item
Suppose $f_{\varphi}$ fixes a subgeometry $[\Sigma_{\varphi}]\cong \PG(3,q)$ of (vector) dimension $4$.
Clearly, $[\Rad(\varphi)]\not\in[{\Sigma}_{\varphi}]$.
  By Lemma~\ref{l2}, $f_{\varphi}$ restricted to $\cH_4=[{\Sigma}_{\varphi}]\cap\cH_5$,
  bijectively maps points of $\cH_4$ into points of $\cH_4$ and fixes $(q^4-1)/(q-1)$
of them. Since every point in the image of $f_{\varphi}$ admits at least $q-1$
preimages in $\cH_5$ distinct from itself
we get
\[ C\geq (q^2-1)(q-1)\mu_4. \]
Plugging this in Equation~\eqref{ewh} and using Corollary~\ref{eqr} and Lemma~\ref{lA2}, we obtain
that for any $q$:
\begin{multline}
\wt(\varphi)\geq \frac{(q^3-q)}{(q^2+1)}\mu_5-\frac{(q^3-q)}{(q^4-1)}A +\frac{q}{(q^4-1)}C\geq\\
 \frac{q(q^2-1)}{(q^2+1)}\mu_5-\frac{(q^3-q)}{(q^4-1)}\xi_5(2) +\frac{q}{(q^4-1)} (q^2-1)(q-1)\mu_4=\\
=\frac{q^{10}-2q^6+q^5-q^4-q^3+2q^2}{q+1}>q^8-q^6.
\end{multline}
\item
  Suppose now that $f_{\varphi}$ does not fix a subgeometry isomorphic to $\PG(3,q)$.
  Then, by Lemma~\ref{v4c1},
$|\Fix(f_{\varphi})|\leq q^2+q+2$.
By Lemma~\ref{tlemma}, $p\in\fA_{\varphi}$ if and only if $f_{\varphi}([p])= [p]$ or
  $f_{\varphi}([p])=0$.

  By Equation~\eqref{eAA}, we have
  \begin{multline*}
    A\leq |\Fix(f_{\varphi})| +(q^2-1)|([\Rad(\varphi)]\cap\cH_{5})|\leq\\
    (q^2+q+2)(q^2-1) + (q^2-1)=
    q^4+q^3+2q^2-q-3.
  \end{multline*}
  We now need to compute a lower bound for $C.$ Since $[\Rad(\varphi)]=\ker(f_{\varphi})$ and we are assuming $\dim(\Rad(\varphi))=1$, the image $\mathrm{Im}(f_{\varphi})$ of the semilinear function $f_{\varphi}$ is a subspace of $\PG(4,q^2)$
  of (vector) dimension  $4$. In particular,
  the image $f_{\varphi}(\cH_5)$ of its restriction to
  $\cH_5$ is a (possibly degenerate) Hermitian surface contained in a projective space $\PG(3,q^2)$.

Let $\cH':=f_{\varphi}(\cH_5)\cap \cH_5.$ By Lemma~\ref{tlemma}, $p\in\fC_{\varphi}$ if
and only if $f_{\varphi}([p])\neq [p]$ and $f_{\varphi}([p])\in\cH'$.
Using the descriptions of intersection
of Hermitian varieties in
\cite{L06,DD07} we see that $|\cH'|\geq q^3+1$.
Thus we get
 \[C\geq (q-1)(q^2-1)(q^3+1).\]
 Plugging this in Equation~\eqref{ewh}, we obtain that for any $q$,
 \begin{multline}
   \label{eeee}
\wt(\varphi)\geq \frac{(q^3-q)}{(q^2+1)}\mu_5-\frac{(q^3-q)}{(q^4-1)}A +\frac{q}{(q^4-1)}C\geq\\
 \frac{q^{10}-q^6+q^5-2q^4-2q^3+2q^2+q}{q^2+1}
\end{multline}
For $q>2$, \eqref{eeee} gives $\wt(\varphi)\geq q^8-q^6$.
And this completes the argument.
For $q=2$ a direct computer search proves that the minimum weight of the code is
once more $192=q^8-q^6$.
\end{enumerate}
\end{itemize}
The above proof directly implies the following characterization of the
minimum weight codewords for $m$ odd.
\begin{corollary}
  \label{codd}
  If either
  \begin{itemize}
    \item $m>5$ is odd or
    \item $m=5$ and $q\not=2$,
    \end{itemize}
    then the minimum weight codewords of $\cC(\HH_{m,2})$
    correspond to bilinear alternating forms $\varphi$ with
    $\dim(\Rad(\varphi))=m-2$ and such that $[\Rad(\varphi)]$ meets
    $\cH_m$ in a Hermitian cone of the form $[\Pi_1]\cH_{m-3}$.
\end{corollary}
\begin{remark}
  For $q=2$ and $m=5$, an exhaustive computer search shows that the
  characterization of Corollary~\ref{codd}
  does not hold, as there are $24948$
  codewords of minimum weight $192$; $5940$ of
  these are associated with bilinear forms with radical
  of dimension $1$ while the remaining $19008$ are  associated with forms
  with radical of dimension $3$. The
  forms with $3$-dimensional radical are as those described in Corollary~\ref{codd}.
  Incidentally, the full list of weights for this code is $0,192,216,224,232,256$.
\end{remark}

\subsection{Minimum distance of $\cH_{m,2}$ with $m$ even}
\label{lev0}
In this section we assume $m$ to be even. Then the Witt index of the Hermitian form $\eta$ is  $n=m/2.$
Let $\varphi$  be an alternating form on $V.$
Recall from Equation~\eqref{ewh} that
\[ \wt(\varphi)=
\frac{q^{2m-7}}{(q^4-1)}((q^2-1)\mu_m - A) +\frac{q^{m-4}}{q^4-1}B.\]
\begin{prop}\label{min wgt m even}
There exists a bilinear alternating form $\varphi$ with $\dim(\Rad(\varphi))=m-2$ such that $\wt(\varphi)=q^{4m-12}$ and $\wt(\varphi')\geq q^{4m-12}$ for any other form $\varphi'$ with $\dim(\Rad(\varphi'))=m-2.$
\end{prop}
\begin{proof}
Let $\varphi$ be a bilinear alternating form with radical of dimension $m-2.$  In order to determine the weight of the word of $\cH_{m,2}$ induced by the form $\varphi$, we need to determine
the number of lines of $\cH_{m,2}$ which are not totally isotropic for $\varphi$.

Since, by hypothesis, the radical of $\varphi$ has dimension $\dim(\Rad(\varphi))=m-2$,  a line $\ell$ is totally isotropic for $\varphi$ if and only if $\ell\cap[\Rad(\varphi)]\neq \emptyset.$
If $S$ denotes the matrix representing $\varphi$, we have $\rank(S)=2$; so,
according to the results obtained in Section~\ref{westim},
we have that for words associated to the value $A_{\max}$ it must be $i=1$ and
$[\Rad(\varphi)]\cap\cH_m=[\Pi_t]\cH_{m-2-t}$ is a degenerate Hermitian variety with radical $[\Pi_t]$ of dimension $t.$
By Equation~\eqref{eq-A} and Lemma~\ref{ll0}, since $m$ is even, the maximum number of points $\mu_{m-2}^{\max}$ of
$[\Rad(\varphi)]\cap \cH_m$ is attained for $t=2i=2$ (in this case $[\Pi_2]$ is a line), i.e. $\mu_{m-2}^{\max}=\mu_{m-2}(2)=q^4\mu_{m-4}+q^2+1$ (last equality comes from Equation~\eqref{numero punti}).
Hence the number of points of $\cH_m\setminus[\Rad(\varphi)]$ (see Equations~\eqref{mu_m} and~\eqref{numero punti}) is at least
\[\mu_m-\mu_{m-2}^{\max}= \mu_{m}-q^4\mu_{m-4}-q^2-1=q^{2m-3}+q^{2m-5}. \]
Assume $[\Rad(\varphi)]\cap\cH_{m}=[\Pi_2]\cH_{m-4}$ and
consider a point $[p]\in\cH_{m}\setminus [\Rad(\varphi)]$; we study
 $[p]^{\perp_{\eta}} \cap \cH_{m}\cap[\Rad(\varphi)]$.
\begin{itemize}
\item
   \framebox{Case $[\Pi_2]\subseteq [p]^{\perp_\eta}$}.
  Note that $[\Rad(\varphi)]=[\Pi_2]^{\perp_\eta}$ because $[\Rad(\varphi)]\cap\cH_m\subseteq [\Pi_2]^{\perp_\eta}$
  and $\dim(\Pi_2^{\perp_\eta})=\dim(\Rad(\varphi))=m-2$ (indeed, if $[x]\in [\Rad(\varphi)]\cap\cH_m$, then
  $[x]\in [\Pi_2]\cH_{m-4}$; hence, $[x]\in [\Pi_2]^{\perp_\eta}$).
   This implies that every point $[p]$ such that $[\Pi_2]\subseteq [p]^{\perp_\eta}$, i.e. $[p]\in [\Pi_2]^{\perp_\eta}$ is also in $[\Rad(\varphi)]$ while we were assuming $[p]\in\cH_{m}\setminus[\Rad(\varphi)].$ Thus, this case can not happen.
\item
  \framebox{Case $[\Pi_2]\not\subseteq [p]^{\perp_\eta}$}.  Then $[p]^{\perp_\eta}\cap(\cH_{m}\cap [\Rad(\varphi)])\cong [\Pi_1]\cH_{m-4}$.
  In this case there are $(\mu_{m-2}-q^2\mu_{m-4}-1)=q^{2m-7}$ lines through $[p]$ disjoint
  from $[\Rad(\varphi)]$.
\end{itemize}
   Since all points of $\cH_{m}$ not in $[\Rad(\varphi)]$ are such that $[\Pi_2]\not\subseteq [p]^{\perp_\eta}$, the total number of lines disjoint from $[\Rad(\varphi)]$ is $\frac{(q^{2m-3}+q^{2m-5})q^{2m-7}}{q^2+1}=q^{4m-12}$, i.e.
   it is always possible to find a bilinear alternating form $\varphi$ with $\dim(\Rad(\varphi))=m-2$ such that $\wt(\varphi)=q^{4m-12}$.
   Observe that
   for  any form with $\dim(\Rad(\varphi))=m-2$ such that $[\Rad(\varphi)]\cap\cH_{m}\not=[\Pi_2]\cH_{m-4}$,
   the number of totally $\eta$-isotropic lines disjoint from $[\Rad(\varphi)]$ is larger than the value obtained
   in the case considered above, so $\wt(\varphi)>q^{4m-12}$.
 \end{proof}
 We claim that $q^{4m-12}$ is actually the minimum weight for $m$ even unless $m=4,6$.  We also compute the minimum weight for $m=4,6$.

As before, let $i=(\rank(S))/2$.
When $i=1$, then by Proposition~\ref{min wgt m even}, the minimum weight  of the codewords   is
$d_1=q^{4m-12}$. Assume $i>1$;  we need to distinguish several cases according to the value of $m$.
\begin{itemize}
\item
 \framebox{$m\geq 12$}. By Corollary~\ref{eqr} and Lemma~\ref{lA2}
  $A_{i}\leq\xi_m(i)\leq\xi_m(2)\leq \xi_m(1).$
  From Case~\ref{AiB}) of Corollary~\ref{eqr} we get
  \[ d_i\geq \frac{q^{2m-7}(\mu_{m}(q^2-1)-\xi_{m}(i))}{q^4-1} \geq \frac{q^{2m-7}(\mu_{m}(q^2-1)-\xi_{m}(2))}{q^4-1}. \]
We will show that
  \[\frac{q^{2m-7}(\mu_{m}(q^2-1)-\xi_{m}(2))}{q^4-1} >
  q^{4m-12}. \]
Actually, by straightforward computations, the above condition becomes
  \[ q^{4m-12}+q^{4m-16}-q^{2m-6}-q^{2m-7}>q^{4m-12}\]
  which is true for all values of $q$.
  \item \framebox{$m=10$}. By Lemma~\ref{lA2} we have that the second largest value of $\xi_{10}(i)$ is for $i=5$:
$\xi_{10}(5)=q^{11}+q^{10}-q-1.$
  By Case~\ref{AiB}) of Corollary~\ref{eqr} we have
\[d_i\geq \frac{q^{13}}{q^2+1}\left(\mu_{10}-\frac{1}{q^2-1}\xi_{10}(5)\right)=
q^{28}+q^{24}-q^{18}-q^{14}>q^{28}.\]
 So, the minimum distance is attained by codewords corresponding to matrices $S$ of rank $2$
and, consequently,  $d_{\min}=q^{4m-12}$.

\item \framebox{$m=8$}. By Lemma~\ref{lA2} we have that the second largest value of $\xi_8(i)$ is for $i=4$: $\xi_8(4)=q^9+q^8-q-1.$
  By Case~\ref{AiB}) of Corollary~\ref{eqr} we have
 \[d_i\geq
 \frac{q^{9}}{q^2+1}\left(\mu_{8}-\frac{1}{q^2-1}\xi_{8}(4)\right)=
   q^{20}+q^{16}-q^{14}-q^{10}>q^{20}.\]
So the minimum distance is attained by codewords corresponding to matrices $S$ of rank $2$
and, consequently,  $d_{\min}=q^{4m-12}=q^{20}$.

\item \framebox{$m=6,4$}. By Lemma~\ref{lA2} we have
 $\xi_4(1)<\xi_4(2)$ and $\xi_6(2)<\xi_6(1)<\xi_6(3)$, hence
 the maximum value of $\xi_m(i)$ is for $i=m/2$, i.e. we have that
 the matrix $S$ has maximum rank $m$ and so it is non-singular.

For $i\neq m/2$, by case~\ref{AiB}) of Corollary~\ref{eqr} we have
\[ d_i\geq d_1>\frac{q^{2m-7}}{q^2+1}\left(\mu_m-\frac{1}{q^2-1}\xi_m(m/2)\right)=
q^{4m-12}-q^{2m-6}=d_{m/2}. \]
We shall show that $q^{4m-12}-q^{2m-6}$ is the actual minimum distance.
\begin{lemma}
\label{lm46}
 If $m=4,6$ then there exists a non-singular alternating form $\varphi$ of $V(m,q^2)$
 such that $|\fA_{\varphi}|=A_{m/2}=(q^m-1)(q+1)$.
\end{lemma}
\begin{proof}
 For any non-singular bilinear alternating form $\varphi$ we have
 $|\fA_{\varphi}|=|\Fix(f_{\varphi})\cap\cH_m|(q^2-1)$, see~\eqref{frakA1-frakA2}.
 Choose $\varphi$ to be a symplectic polarity which permutes with $\eta$ (i.e.
 $[u]^{\perp_{\varphi}\perp_{\eta}}=[u]^{\perp_{\eta}\perp_{\varphi}}$ for all
 $[u]\in\PG(m-1,q^2)$). Then $f_{\varphi}(\cH_m)=\cH_m$ and by~\cite[\S 74]{S65},
 $\Fix(f_{\varphi})\cong\PG(m-1,q)$ is a subgeometry over $\FF_q$
  fully contained in $\cH_m$. Hence,
  \begin{multline*} |\fA_{\varphi}|=|\Fix(f_{\varphi})\cap\cH_m|(q^2-1)=|\Fix(f_{\varphi})|(q^2-1)=\\
  \frac{q^m-1}{q-1}(q^2-1)=(q^m-1)(q+1).
\end{multline*}
\end{proof}
By Lemma~\ref{l2}, the bilinear alternating form $\varphi$ given by Lemma~\ref{lm46} is such that
$|\fB_{\varphi}|=0$. Hence, since $A_{m/2}>A_i$ for all $i\neq m/2$ and $A_{m/2}=(q^m-1)(q+1)$, we have
$\wt(\varphi)=q^{4m-12}-q^{2m-6}$. By Equation~\eqref{e5}, as $\wt(\varphi)=q^{4m-12}-q^{2m-6}$, it
follows that $d_{\min}=q^{4m-12}-q^{2m-6}$.
 \end{itemize}
 By the arguments presented before we have the following characterization of
 the minimum weight codewords for $m$ even.
 \begin{corollary}
   \label{ceven}
   If $m=4$ or $m=6$, then the minimum weight codewords of $\cC(\HH_{m,2})$
   correspond to bilinear alternating forms $\varphi$ which are permutable
   with the given Hermitian form $\eta.$
   If $m>6$ is even,
   then the minimum weight codewords
   correspond to bilinear alternating forms $\varphi$ with
   $\dim(\Rad(\varphi))=m-2$ and such that $[\Rad(\varphi)]$ meets
   $\cH_m$ in a Hermitian cone of the form $[\Pi_2]\cH_{m-4}$.
\end{corollary}

Sections~\ref{odd_m} and \ref{lev0} complete the proof of the Main Theorem.

\section*{Acknowledgments}
Both authors are affiliated with GNSAGA of INdAM (Italy) whose support they
acknowledge.

\end{document}


\customlabel{tlemma}{3.2}
\customlabel{l2}{3.3}
\customlabel{eww}{13}
\customlabel{eqr}{3.6}
\customlabel{lA2}{3.7}
\customlabel{min wgt m odd}{3.9}
\customlabel{f_varphi}{15}
\customlabel{eAA}{19}
\customlabel{codd}{3.11}
\customlabel{v4c1}{3.10}
\begin{frontmatter}
  \title{Corrigendum to ``Line Hermitian Grassmann Codes and their Parameters''}
\author[IC]{Ilaria Cardinali}
\ead{ilaria.cardinali@unisi.it}
\address[IC]{Department of Information Engineering and Mathematics, University of Siena,
Via Roma 56, I-53100, Siena, Italy}
\author[LG]{Luca Giuzzi\corref{cor1}}
\ead{luca.giuzzi@unibs.it}
\address[LG]{DICATAM - Section of Mathematics,
  Universit\`a di Brescia,
Via Branze 43, I-25123, Brescia, Italy}
\cortext[cor1]{Corresponding author.}
\begin{abstract}
  In this note we provide a correction to the statement of the main
  result of~\cite{LHG}
  for the case $m=5$.
\end{abstract}
\begin{keyword}
  Hermitian variety \sep Polar Grassmannian \sep Projective Code.
  \MSC[2010] 14M15 \sep 94B27 \sep 94B05
\end{keyword}
\end{frontmatter}
\section{Introduction}
In~\cite{LHG} on page 419, it is wrongly stated that the fixed points
of a properly semilinear collineation $\varphi$ of $\PG(m-1,q^2)$ must be contained in a subgeometry $\Sigma_{\varphi}\cong\PG(t,q)$ for some $t$ with $0\leq t\leq m-1$.
This assertion is not correct; what can be shown is 
that the maximum number of fixed points of $\varphi$ is at most
$\frac{q^t-1}{q-1}=|\PG(t,q)|$ for some $t$ with $0\leq t\leq m-1$; see
Theorem~\ref{new} below for the precise statement.
This mistake does not affect the statement of the main result of the paper,
i.e. that all of the codes under investigation have minimum
distance $d_{\min}=q^{4m-12}-q^{3m-9}$. The proof of Lemma~\ref{l2} is
minimally affected by the mistake, as the assertion
\emph{``and the semilinear transformation $f_{\varphi}$ fixes a subgeometry [\dots] $\PG(m-1,q)$ of maximal dimension''} should be replaced with
\emph{``and, by Theorem~\ref{new}, the fixed points of $f_{\varphi}$ span $\PG(V)$''}.
With this provision, all of the arguments for the cases $m\neq5$
remain exactly the same.

In the case $m=5$ however, the proof that the minimum distance is $q^8-q^6$
must be amended. Also, it can be shown that for $m=5$ and any value of $q$
(not just for $q=2$ as stated in Corollary~\ref{codd})
there are minimum weight codewords associated with alternating forms
with radical of dimension $1$; see Remark~\ref{r3}.
We also point out that in Corollary~\ref{codd} the number of words of minimum weight
with radical of dimension $1$ and those with radical of dimension $3$
has been swapped and the correct statement is Remark~\ref{r4}.

Throughout this corrigendum we shall use the notation and theorem numbers
of~\cite{LHG}. The new statements are labeled with roman letters in order
to avoid confusion. 

In Section~\ref{sec2} we recall and prove
some statements on the number and geometry of the fixed
points of a semilinear collineation; in particular Theorem~\ref{new}
should replace the unnumbered statement on page 419; Lemma~\ref{v4c1new}
is the weaker version of Lemma~\ref{v4c1} which shall be used to prove the main result.

In Section~\ref{sec3} we reinstate the main theorem of the paper in
the special case $m=5$ as
Theorem~\ref{mth} and provide a corrected proof.
We conclude with some remarks about minimum weight codewords for $m=5$.

The notation is the same as of~\cite{LHG}.

\section{Fixed points of a semilinear collineation}
\label{sec2}

Throughout the paper we need to estimate the maximum number of fixed points
of a semilinear collineation of $\PG(m-1,q^2)$. In the original version of the paper it was
stated (without proof)
that these points must belong to some subgeometry; so we would
immediately have that they are at most $\frac{q^{m}-1}{q-1}$. As mentioned
in the introduction, this is not true in general.
The correct statement is the 
following, whose proof comes from~\cite[Theorem 3.8]{parte2} and~\cite[Theorem 3.9]{parte2}.

  \begin{theorem}
\label{new}
  Let $f_{\varphi}$ be a (possibly degenerate) proper semilinear transformation
  of $\PG(m-1,q^2)\to\PG(m-1,q^2)$. Then,
  \begin{enumerate}
    \item the fixed points of $f_{\varphi}$
      are contained in a subspace $\Sigma_{\varphi}$ which meets
      $[\Rad(\varphi)]$ trivally.
    \item Put $r=\dim(\Sigma_{\varphi})$. 
      The number of fixed points of $f_{\varphi}$ is at most
      $\frac{q^r-1}{q-1}$.
    \end{enumerate}
  \end{theorem}

 Let $V,W$ be vector spaces over the same field;
 throughout this section we shall denote by the same symbol
 a semilinear collineation $f:\PG(V)\to\PG(W)$ and
 the corresponding semilinear map $f:V\to W$.
 
 \begin{lemma}  \label{condizione (*)}
   Let $f_{\varphi}\colon \PG(4,q^2)\to\PG(4,q^2)$ be
   a semilinear collineation.
   Suppose $[v_1]$,
   $[v_2]$, $[v_3]$, $[v_4]$ are four indipendent points fixed by $f_{\varphi}$. Then the following hold.
   \begin{enumerate}
     \item There exist $\alpha_i\in\FF_{q^2}^*$ for $i=1,\dots,4$ such
       that $f_{\varphi}(v_i)=\alpha_i v_i$;
     \item
       Take $K\subseteq\{1,2,3,4\}$ with $K\neq\emptyset$. There exists a fixed point $[c]:=[\sum_{i=1}^4c_iv_i]\in [\langle v_1,v_2,v_3,v_4\rangle]$
       with $c_k\neq0$ $\forall k\in K$
       if and only if
       \begin{equation}
         \label{cond30}\alpha_i^{q+1}=\alpha_j^{q+1},\quad
         \forall i,j\in K.
       \end{equation}
     \end{enumerate}
     \end{lemma}
     \begin{proof}
       Since the points $[v_i]$ are fixed by $f_{\varphi}$
       we have
       $f_{\varphi}([v_i])=[v_i]$, whence
       $f_{\varphi}(v_i)=\alpha_i v_i$, with $\alpha_i\in\FF_{q^2}^*$ for
       all $i$. This proves point 1.

       Let $\Pi:=[\langle v_1,v_2,v_3,v_4\rangle]$.
       Then $f_{\varphi}$ restricted to $\Pi$ is invertible by Theorem~\ref{new}.
       Take $[c]=[\sum_{i=1}^nc_iv_i]$ be a point fixed by $f_{\varphi}$
       in $\Pi$.
       Then
       \[ [\sum_{i=1}^r c_i^q\alpha_i v_i]= [\sum_{i=1}^4c_i^qf_{\varphi}(v_i)]=
          f_{\varphi}([c])=
         [c]=
         [\sum_{i=1}^4 c_iv_i].
       \]
       So, there is $\lambda\in\FF_{q^2}^*$ such that
       \begin{equation}
         \label{c0} c_i^q\alpha_i=\lambda c_i,\quad \forall i=1,\dots,4.
       \end{equation}

       Suppose first $c_k\neq0$ for all $k\in K$; then,
       from Equation~\eqref{c0} we get
       \[\lambda=\alpha_kc_{i}^{q-1},\quad \forall k\in K; \]
       so, in particular we have
$\alpha_ic_i^{q-1}=\alpha_jc_j^{q-1}$ for all  $i,j\in K$.
Thus, rising to the $(q+1)$-th power,
\[
\alpha_i^{q+1}=\alpha_j^{q+1},\quad \forall i,j\in K.
\]

Conversely, suppose that~\eqref{cond30} holds. Since $K\not=\emptyset$, we can assume without
loss of generality $1\in K$. Hence, for all $i\in K$ we have
$(\alpha_i^{-1}\alpha_1^{1})^{q+1}=1$, whence there exist some scalars $c_i\in \FF_{q^2}^*$ such that $\alpha_i^{-1}\alpha_1=c_i^{q-1}.$
Put $c:=\sum_{i\in K} c_iv_i$. We claim that $[c]$ is fixed by
$f_{\varphi}$. Indeed, we have
\[ f_{\varphi}([c])=f_{\varphi}([\sum_{i\in K}c_iv_i])=[\sum_{i\in K}c_i^q\alpha_i
  v_i].
  \]
  Observe that $f_{\varphi}([c])=[c]$ if and only if there exists $\lambda\neq0$ such  that $c_i^q\alpha_i=\lambda c_i$ for all $i\in K$.
  However, since $c_i\neq0$ for all $i\in K$, the above condition
  is the same as to say that there exists $\lambda\in\FF_{q^2}^*$ such
  that $c_i^{q-1}\alpha_i=\lambda$. By definition of $c_i$,
  this equation becomes $\alpha_1=(\alpha_i^{-1}\alpha_1)\alpha_i=\lambda.$ So, taking $\lambda=\alpha_1$
   the thesis follows.
\end{proof}

With the same assumption of Lemma~\ref{condizione (*)},  we have
\begin{lemma}\label{condizione Delta}
 Suppose that there exists
  a $T\subseteq\{1,2,3,4\}$ and $\gamma\in\FF_q^*$ such that
  $\alpha_i^{q+1}=\gamma$ for any $i\in T$ and there exists
  $r\in\{1,2,3,4\}\setminus T$ with $\alpha_r^{q+1}\neq\gamma.$
 Let $[c]=[\sum_{i=1}^4c_iv_i]$ be a fixed point for $f_{\varphi}$.

  Under these assumptions,
  if $c_r\neq0$ then $c_i=0, \,\forall i\in T.$
\end{lemma}
\begin{proof}
Assume by way of contradiction that
$\exists i\in T\colon c_i\neq 0\neq c_r$. By Lemma~\ref{condizione (*)} this
implies that $\{i,r\}\subseteq\{1,2,3,4\}$ and
$\alpha_i^{q+1}=\alpha_r^{q+1}$, against the hypothesis
$\gamma=\alpha_i^{q+1}\neq\alpha_r^{q+1}$.
It follows that $c_i=0$ for all $i\in T$.
\end{proof}

In light of the above correction,
the following must replace~\cite[Lemma~\ref{v4c1}]{LHG}.

\begin{lemma}
  \label{v4c1new}
  Let $f_{\varphi}\colon \PG(4,q^2)\to\PG(4,q^2)$ be the semilinear collineation $f_{\varphi}([x])=[x]^{\perp_{\varphi\perp_{\eta}}} $ given by~\eqref{f_varphi} and let $[\Fix(f_{\varphi})]$ be the set of points fixed by $f_{\varphi}$.
  Then $[\Fix(f_{\varphi})]$ is contained in a  $3$-dimensional projective space $\PG(3,q^2)$ and either $|[\Fix(f_{\varphi})]|=q^3+q^2+q+1$ or
  $|[\Fix(f_{\varphi})]|\leq q^2+q+2$.

  Moreover, if $|[\Fix(f_{\varphi})]|<q^2+q+1$, then
    $|[\Fix(f_{\varphi})]|\leq 2q+2$.
\end{lemma}
\begin{proof}
  By construction, the semilinear collineation $f_{\varphi}$
  has a kernel of odd vector dimension, corresponding to the radical
  of the alternating bilinear form $\varphi$.
  Let $\Pi:=\langle\Fix(f_{\varphi})\rangle$ be the span of the
  vector representatives of the
  points fixed by
  $f_\varphi$.
  By Theorem~\ref{new}, $\ker(f_{\varphi})$ meets
  $\Pi$ trivially. In particular, if
  $\dim(\Rad(\varphi))=3$ then, by Theorem~\ref{new}, $\dim(\Pi)\leq2$
  and, again by the same theorem,
  $|[\Fix(f_{\varphi})]|\leq q+1$. Clearly, in this case $[\Pi]$ is
  contained in some $3$-dimensional projective space $\PG(3,q^2)$.

  If  $\dim(\Rad(\varphi))=1$, again
  by Theorem~\ref{new}, $\dim(\Pi)\leq 4$, i.e.
  $[\Pi]$ is contained in a projective space $\PG(3,q^2)$.
  This proves the first part of the thesis.

  Put $g:=f_{\varphi}|_{[\Pi]}$. By
  construction, $g$ is non degenerate and
  $\Fix(f_{\varphi})=\Fix(g)$.


  Suppose  that $\dim(\Pi)=4$.
  By Theorem~\ref{new} we have
  $|[\Fix(g)]|\leq (q^4-1)/(q-1)=q^3+q^2+q+1$.
  Fix now a basis of $\Pi$ consisting of vectors $v_1,v_2,v_3,v_4$ corresponding
  to $4$ independent fixed points of $g$.
   Then, there exist
  $\alpha_i\in \FF_{q^2}^*$ such that $g(v_i)=\alpha_i v_i$ for $i=1,\dots,4$.

  Let $\Omega:=\{ \alpha^{q+1}_i : i\in 1,\dots,4 \}$.

 Clearly, $1\leq |\Omega|\leq 4$.

  Let now $[c]=[\sum_{i=1}^4c_iv_i]$ be a point fixed by $g$.
  The components $c_i$ of $c$ must satisfy~\eqref{c0},
  that is, there exists $\lambda\in\FF_{q^2}^*$ such that
  the following system has solution:
  \begin{equation}
    \label{ee0a}
    \begin{cases}
    \alpha_1 c_1^q= \lambda c_1 \\
    \alpha_2 c_2^q=\lambda c_2 \\
    \alpha_3 c_3^q=\lambda c_3 \\
    \alpha_4 c_4^q=\lambda c_4 \\
  \end{cases}
\end{equation}
We now distinguish $4$ cases according to the size of $\Omega$.
 \begin{enumerate}[(a)]
 \item \fbox{$|\Omega|=1$.} \\ In this case we
   can take $K=\{1,2,3,4\}$ in Lemma \ref{condizione (*)}. So,
   there exists  a fixed point $[c]$ such that $0\not\in\{c_1,c_2,c_3,c_4\}$.
    Note that if $(c_1,c_2,c_3,c_4)$ is a solution to~\eqref{ee0a},
    then for any $(\omega_1,\dots,\omega_4)\in\FF_q^4$
    also $(\omega_1 c_1,\omega_2 c_2,\omega_3 c_3,\omega_4 c_4)$ is a solution.
      So the number of non-proportional, non-zero solutions to~\eqref{ee0a}
   (i.e. distinct point representatives) is at least $(q^4-1)/(q-1)=
   q^3+q^2+q+1$. This is also the maximum number of fixed points of $g$.
   It follows that $|[\Fix(g)]|=q^3+q^2+q+1.$

 \item\label{Om2} \fbox{$|\Omega|=2.$} \\
   Without loss of generality, up to permuting the indexes of the vectors, we
   have two subcases:
   \begin{enumerate}[{b}.i)]
   \item\label{Om2c1}
     $\alpha_1^{q+1}=\alpha_2^{q+1}=\alpha_3^{q+1}\neq\alpha_4^{q+1}$.
     In this case, put $T=\{1,2,3\}$ and take $r=4$.
     Then, by Lemma~\ref{condizione Delta} we have that either
     $(c_1,c_2,c_3)=(0,0,0)$ (and $c_4\not=0$)or $c_4=0$.
     In particular, either
     $[c]$ is contained in the subspace spanned by $[v_1],[v_2],[v_3]$
     or $[c]=[c_4]$.
     By Theorem~\ref{new} it follows that $|[\Fix(g)]|\leq q^2+q+1+1=
     q^2+q+2$.
   \item
     Suppose now $\alpha_1^{q+1}=\alpha_2^{q+1}\neq\alpha_3^{q+1}=\alpha_4^{q+1}$.
     Since $[c]$ is fixed take $T=\{1,2\}$ and $r\in\{3,4\}$ in
     Lemma~\ref{condizione Delta}.
     It follows that if $c_3\neq 0$ or $c_4\neq0$, then $(c_1,c_2)=(0,0)$.
     Similarly, taking $T=\{3,4\}$ and $r\in\{1,2\}$ in
     Lemma~\ref{condizione Delta} we have that
     if $c_1\neq0$ or $c_2\neq0$, then $(c_3,c_4)=(0,0)$.
     In particular   $[c]$ belongs to
     the union $[\langle v_1,v_2\rangle]\cup[\langle v_3,v_4\rangle]$.
     By Theorem~\ref{new}, each of the lines $[\langle v_1,v_2\rangle]$
     and $[\langle v_3,v_4\rangle]$ contains at most $q+1$ fixed points.
     So $|[\Fix(g)]|\leq 2(q+1)$.
\end{enumerate}
     \item \fbox{$|\Omega|=3.$} \\
       Without loss of generality we can suppose $\alpha_1^{q+1}=\alpha_2^{q+1}$
       and $\alpha_3^{q+1}\neq\alpha_4^{q+1}$ with $\alpha_4^{q+1}\neq
       \alpha_1^{q+1}\not= \alpha_3^{q+1}$.
       Using Lemma~\ref{condizione Delta} with $T=\{1,2\}$ and $r\in\{3,4\}$ we see that
       if $c_3\neq0$ or $c_4\neq0$ then $(c_1,c_2)=(0,0)$. This is
       equivalent to say that
       $c_1\neq0$ or $c_2\neq0$ implies $(c_3,c_4)=(0,0)$.
       So, we have $[c]\in[\langle v_1,v_2\rangle]\cup[\langle v_3,v_4\rangle]$.
       Suppose now $[c]\in[\langle v_3,v_4\rangle]$.
       If there were $[c]=[c_3v_3+c_4v_4]$ with $c_3\neq0\neq c_4$ fixed we would
       have $\alpha_3=\lambda c_3^{1-q}$ and $\alpha_4=\lambda c_4^{1-q}$,
       whence $\alpha_3^{q+1}=\alpha_4^{q+1}$, contradicting the conditions
       on the $\alpha_i$'s. So either $c_3=0$ or $c_4=0$, that is
       $[c]=[v_3]$ or $[c]=[v_4]$.  As $[v_3]$ and $[v_4]$ are fixed by $g$,
       there are exactly two fixed points for $g$ in $[\langle v_3,v_4\rangle]$.
       By Theorem~\ref{new}, there are at most $q+1$ fixed points for $g$
       on the line $[\langle v_1,v_2\rangle]$.
       It follows that
       $|[\Fix(g)]|\leq (q+1)+2=q+3$.
     \item \fbox{$|\Omega|=4.$} \\
    Suppose in this case $i,j\in\{1,\dots,4\}$ and $i\neq j$.
     If there were a point $[c]$ with both $c_i\neq0$ and $c_j\neq0$,
     then
     \[ \alpha_i=\lambda c_i^{1-q},\qquad \alpha_j=\lambda c_j^{1-q}
       ,\]
     whence $\alpha_i^{q+1}=\alpha_j^{q+1}$, which is a contradiction.
     It follows that exactly one of
     the coordinates $c_i$ of $c$ is different from $0$.
     So $[c]\in\{ [v_1], [v_2], [v_3], [v_4] \}$ and the number of
     fixed points is $4$.
   \end{enumerate}
  If the vector dimension of
  $\Pi$ is strictly less than $4$, then, by Theorem~\ref{new},
  the number of fixed points of $g$ is at most $(q^3-1)/(q-1)=q^2+q+1$.

Performing now the same
  analysis just done for  $\dim(\Pi)=4$ but taking as $\Omega$  the set associated to three
  fixed points whose vector representatives span $\Pi$,
  it is not difficult to see that if $\dim(\Pi)=3$, we have
  that the only possible numbers of fixed points are $q^2+q+1$,
  $q+2$ and $3$.

  Analogously, if $\dim(\Pi)=2$, then the number of fixed points
  is either $q+1$ or $2$ and  finally if $\dim(\Pi)=1$ there is just
  a single fixed point.

     By comparing all the obtained values for $|\Fix(f_{\varphi})|$,  the thesis follows.
\end{proof}

\begin{lemma}\label{nuovo lemma}
  Let $[\Pi]$ be a plane of $\PG(4,q^2)$ spanned by some
  fixed points of $f_{\varphi}$.
  Put $g:=f_{\varphi}|_{[\Pi]}$.
  Then the following hold.
  \begin{enumerate}
  \item
    $g\colon [\Pi]\rightarrow [\Pi], g(x)=x^{\perp_{\varphi}\perp_{\eta}}$
    is a non-degenerate semilinear map.
  \item Either $|[\Fix(g)]|=q^2+q+1$ or $|[\Fix(g)]|\leq q+2$.
  \item If $|\Fix(g)]|=q^2+q+1$,
    then $[\Fix(g)]$ is a Baer
    subgeometry ${\mathfrak G}\cong\PG(2,q)$ contained in
    $\PG(2,q^2)$.
  \end{enumerate}
\end{lemma}
\begin{proof}
  \begin{enumerate}
  \item The first part of the lemma follows from Theorem~\ref{new},
    since the plane $[\Pi]$ is contained in the space spanned
    by the fixed points of $f_{\varphi}$ and $g([x])=f([x])$ for all
    $[x]\in[\Pi]$.
  \item This follows from the same arguments used in the proof of
     Lemma~\ref{v4c1new}.
  \item
    If $g$ fixes
    $q^2+q+1$ points in $\PG(2,q^2)$, then these points
    cannot be contained in a line; so
    $g$
    must fix
    at least three independent
    points $[v_1]$, $[v_2]$, $[v_3]$ and we have
    \[ g(v_1)=\alpha_1 v_1,\quad g(v_2)=\alpha_2 v_2,
      \quad g(v_3)=\alpha_1 v_3. \]
    Using
    the same arguments as in the proof of Lemma~\ref{v4c1new},
    we see that
    $\alpha_1^{q+1}=\alpha_2^{q+1}=\alpha_3^{q+1}$.
      Now, for $i=1,2,3$, choose $c_i\in\FF_{q^2}^*$ such that
      $c_i^{q-1}=(\alpha_1\alpha_i^{-1})$.

      Then, for any
      $(\delta_1,\delta_2,\delta_3)\in\FF_q^3$, the point
      $[c_{(\delta_1,\delta_2,\delta_3)}]:=
      [\delta_1c_1v_1+\delta_2c_2v_2+\delta_3c_3v_3]$
      is fixed by $g$. We have
      $[c_{(\delta_1,\delta_2,\delta_3)}]=[c_{(\delta_1',\delta_2',\delta_3')}]$
      if and only if
      $(\delta_1,\delta_2,\delta_3)=\lambda
      (\delta_1',\delta_2',\delta_3')$
      for some non-zero $\lambda\in\FF_q^*$.

      In this way we obtain all
      $q^2+q+1$ distinct fixed points of $g$, so this gives
      all of the fixed points of $g$ and there
      is a bijective embedding
      \[ \PG(2,q)\cong\frac{\FF_q^3\setminus\{0\}}{\FF_q^*}\to
      [\Fix(g)]\subseteq \PG(2,q^2). \]
    It follows that $[\Fix(g)]$ is a subplane of $\PG(2,q^2)$.
    \end{enumerate}
\end{proof}
\begin{corollary}
\label{cofix}
  Suppose $|[\Fix(f_{\varphi})]|=q^2+q+2$. Then
  $[\Fix(f_{\varphi})]={\mathfrak G}\cup\{[v_4]\}$
  where ${\mathfrak G}$ is a Baer subplane contained in
  a plane $[\Pi]$ and $[v_4]$ is a point not in $[\Pi]$.
\end{corollary}
\begin{proof}
  Since $|[\Fix(f_{\varphi})]|>q^2+q+1$, we have
  $\dim\langle\Fix(f_{\varphi})\rangle=4$.
  We now adopt the notation of the proof of Lemma~\ref{v4c1new}. Accordingly, $v_1,\,v_2,\,v_3,\,v_4$ are four linearly independent vectors corresponding to points fixed by $f_{\varphi}.$

  By the analysis performed in the proof of Lemma~\ref{v4c1new},
  we are in case \ref{Om2}.i) of point (\ref{Om2}), so
  there are $q^2+q+1$ fixed points of $f_{\varphi}$ contained
  in the plane $[\Pi]$ spanned by  $[v_1],[v_2]$ and $[v_3]$ and the futher point $[v_4]\not\in \Pi.$

  By Lemma~\ref{nuovo lemma} these $q^2+q+1$ fixed points are
  a Baer subgeometry $\mathfrak G$.
  As $[v_4]$ is also a fixed point and $[v_4]\not\in[\Pi]$
  we have $[\Fix(f_{\varphi})]=\mathfrak{G}\cup\{[v_4]\}$.
\end{proof}

\section{Main theorem for $m=5$}
\label{sec3}
In this section we reinstate the main theorem of the paper for the
case $m=5$ and provide its corrected proof. 

We need a further preliminary lemma of geometric nature.
\begin{lemma}
  \label{IntH4}
  Let $\cH_4$ and $\cH_4'$ be two Hermitian surfaces in $\PG(3,q^2)$
  not containing planes
  and suppose that $\cH_4\cap\cH_4'$ contains a degenerate
  Hermitian curve $\Lambda$, i.e. a set of $q+1$
  lines through a point $[v]$ contained in a plane $\Sigma$.
  Then, either
  $\cH_4\cap\cH_4'=\Lambda$
  or $|\cH_4\cap\cH_4'|\geq(q^2-1)^2$.
\end{lemma}
\begin{proof}
  If $\cH_4\cap\cH_4=\Lambda$, then there is nothing to prove.
  So suppose $\cH_4\cap\cH_4'\neq\Lambda$.

  If the linear system determined by the Hermitian surfaces
  $\cH_4$ and $\cH_4'$ contains at least
  one non-degenerate Hermitian variety,
  since $\cH_4\cap\cH_4'\neq\Lambda$, then $\cH_4\cap\cH_4'$ is
  not contained in a plane; so, from~\cite{DD07} or \cite{L06}
  we get $|\cH_4\cap\cH_4'|\geq(q^4-1)$.

  Suppose now that all the Hermitian varieties in the linear system
  generated by $\cH_4$ and $\cH_4'$ are degenerate and
  $\Lambda\subsetneq\cH_4\cap\cH_4'$.
  As $\cH_4$ and $\cH_4'$ do not contain planes, we can assume
  that $\cH_4$ and $\cH_4'$ are two
  Hermitian cones with vertices respectively
  points $[v]$ and $[v']$ and basis suitable non-degenerate
  Hermitian curves contained in a common plane $\Theta$
  not through $[v]$ and $[v']$.
  As all lines contained in a Hermitian cone pass through its vertex
  and $\cH_4\cap\cH_4'$ contains at least $q+1>1$ lines, it must
  necessarily be $[v]=[v']$.
  It follows that $\cH_4\cap\cH_4'$ is a cone of vertex $[v]$
  projecting the intersection of two non-degenerate Hermitian
  curves lying in $\Theta\neq\Sigma$.

  Since $\cH_4\cap\cH_4'$ is a cone, for any point $[p]\neq [v]$
  with $[p]\in\cH_4\cap\cH_4'$, we have
  $[\langle v,p\rangle]\subseteq\cH_4\cap\cH_4'$.
  In particular, if there is $[p]\in(\cH\cap\cH')\setminus\Lambda$,
  then there is also at least one point
  $[p']\in\Theta\cap(\cH_4\cap\cH_4')\setminus\Lambda$.
  Since $\Theta\cap\Lambda$ is a subline $\ell'$ consisting of $(q+1)$
  points, we have that in this case
  $\Theta\cap\cH_4\cap\cH_4'$
  must contain $\ell'\cup\{[p']\}$.
  Looking at~\cite{Kestenband81} we see that this
  condition rules out cases (c) in which the intersection
  is just a subline, (f) where the intersection is
  a single point and (g) in which the intersection is
  an arc. By analyzing the remaining cases,
  $\Theta\cap\cH_4\cap\cH_4'$ must contain at least
  $(q^2+1)$ points. Consequently,
  $|\cH_4\cap\cH_4'|\geq q^2(q^2+1)+1=q^4+q^2+1>(q^4-1)^2$.
  This completes the proof.
\end{proof}

\begin{theorem}[Main theorem for $m=5$]
    \label{mth}
    A line Hermitian Grassmann code deﬁned by a non-degenerate Hermitian form
    on a vector space $V(5,q^2)$ is a $[N,K,d_{\min}]$ code where
    \[ N=q^{10}-q^8+q^7-q^3+q^2-1,\qquad K=10, \]
    \[ d_{\min}=q^8-q^6. \]
\end{theorem}

\begin{proof}
  By Proposition~\ref{min wgt m odd}, there exists an alternating bilinear form $\varphi$ with $\dim(\Rad(\varphi))=3$ such that $\wt(\varphi)=q^8-q^6$ and any other form $\varphi$ with $\dim(\Rad(\varphi))=3$ has weight greater then $q^8-q^6.$ So, we need to show that there
 are no alternating bilinear forms with $\dim(\Rad(\varphi))=1$ having weight
 less than $q^8-q^6$.
 Assume henceforth that $\varphi$ is an alternating bilinear form with $\dim(\Rad(\varphi))=1$. We shall determine a lower bound $d_2$ for the weights $\wt(\varphi)$ and prove that  $\wt(\varphi)\geq d_2\geq q^8-q^6.$

 We need to estimate the values $A,B,C$ in order to apply equation~\eqref{eww} for $m=5$, i.e.
\begin{equation}\label{nuova pesi m=5}
  \wt(\varphi)=\frac{(q^{3}-q)}{(q^2+1)}\mu_{5} - \frac{(q^{3}-q)}{(q^4-1)}A + \frac{q}{(q^4-1)}C.
\end{equation}

 Recall that $C=|\fC_{\varphi}|.$
By Lemma~\ref{tlemma}, $p\in\fC_{\varphi}$ if and only if $f_{\varphi}([p])\neq [p]$ and
$f_{\varphi}([p])\in\cH_5$.
If $[p]=[\Rad(\varphi)]$, then $p\in\fA_{\varphi}$ and this is not possible because $\fA_{\varphi}$ and $\fC_{\varphi}$ are disjoint.
So, suppose $[p]\neq [\Rad(\varphi)]$.
Since $[\Rad(\varphi)]=\ker(f_{\varphi})$, if $[p]\in\cH_5$ and $f_{\varphi}([p])=[x]\in\cH_5$, then $f_{\varphi}([p+\alpha \Rad({\varphi})])=[x]\in\cH_5$ for
any $\alpha\in\FF_{q^2}$.
On the other hand, for any given point $[p]\in\cH_5$, the line $[p,\Rad(\varphi)]$ meets $\cH_5$ in either $1$, $(q+1)$ or
$(q^2+1)$ points. Also if $[p,\Rad(\varphi)]$ meets $\cH_5$ in $1$ or
$(q^2+1)$ points, then $[p]\in[\Rad(\varphi)]^{\perp_{\eta}}$.

In particular, take
$[x]\in(\cH_5\cap f_{\varphi}(\cH_5))$ and let $[p]$ be a preimage of $[x]$
in $\cH_5$
(possibly $[x]$ itself if $[x]$ is fixed). Then, if $[p]\in\cH_5\setminus
[\Rad(\varphi)]^{\perp_{\eta}}$, we have that $[x]$ admits at least
$q-1$ preimages distinct from $[x]$.
If instead $[p]\in\cH_5\cap [\Rad(\varphi)]^{\perp_{\eta}}$ then
either $[x]$ has at most one preimage in $\cH_5$ which is $[p]$,
in case $[\Rad(\varphi)]\not\in\cH_5$ or $[x]$ has at least
$q^2-1$ preimages in $\cH_5$ distinct from itself if $[\Rad(\varphi)]\in\cH_5$.

In light of Lemma~\ref{v4c1new},
observe that
\begin{itemize}
\item
  $|[\Fix(f_{\varphi})]|<q^3+q^2+q+1$ implies $|[\Fix(f_{\varphi})]|\leq q^2+q+2$ and
\item
  $|[\Fix(f_{\varphi})]|<q^2+q+1$ implies $|[\Fix(f_{\varphi})]|\leq 2q+2$.
\end{itemize}

We now distinguish three cases, depending on the number
of fixed points of $f_{\varphi}$:
\begin{enumerate}[a)]
\item\label{case-a}
  \fbox{$|[\Fix(f_{\varphi})]|=q^3+q^2+q+1$}.\\
  Since $\dim(\Rad({\varphi}))=1$, put $[R]:=[\Rad(\varphi)].$
  If $[R]\not\in \cH_5$, then  $A\leq \frac{(q^4-1)(q^2-1)}{q-1}$ (since there could exist points fixed by $f_{\varphi}$ and not contained in $\cH_5$).
  If $[R]\in\cH_5$, we have  $A\leq \frac{(q^4-1)(q^2-1)}{q-1}+(q^2-1)$.
  By Theorem~\ref{new}, there exists a subspace $[\Pi]\cong \PG(3,q^2)$ disjoint from $[R]$ and containing $[\Fix(f_{\varphi})].$
Since $f_{\varphi}$ restricted to $[\Pi]$ is non-singular and it fixes $(q^4-1)/(q-1)>(q^3-1)/(q-1)$  points, it follows from~Theorem~\ref{new} that $[\Pi]=[\langle  \Fix(f_{\varphi})\rangle]$.

Let ${\mathbb B}:=([v_1],\dots,[v_4])$ be a basis of $\Pi$ composed of
vector representatives of fixed points for $f_{\varphi}|_{[\Pi]}$.
  Put $g:=f_{\varphi}|_{[\Pi]}.$
  Hence $g(v_i)=\alpha_i v_i$ for some scalars $\alpha_i\in \FF_{q^2}\setminus \{0\}$.
   By the cases (b), (c), (d) of the  proof of Lemma~\ref{v4c1new}, if $|\{\alpha_i^{q+1}\}_{i=1}^4|>1$, then the number of points fixed by $g$ would be smaller then $(q^4-1)/(q-1)$.
   Hence there exists $\gamma\in \FF_{q^2}\setminus \{0\}$ such that $\gamma=\alpha_i^{q+1}$
  for all $i=1,\dots, 4$.
  By Lemma~\ref{condizione (*)}, there exists a point $[c]$ fixed by $g$ in $[\Pi]$ with all non-zero
  components with respect to ${\mathbb B}$.

  As the fixed points of $g$ are also fixed points of
  $g^2$, we have that
  $g^2$ is a linear transformation (recall that the automorphism $\varphi$ is involutory) fixing at least
  a frame of $[\Pi]$ consisting of ${\mathbb B}\cup\{[c]\}$,
  hence $g^2=\mathit{Id}$. This implies that $g$ is an involutory semilinear collineation.

  Since $g$ is involutory, we have
  \[ g([x])=[x]^{\perp_{\varphi}\perp_{\eta}}=(g^{-1}([x]))=
    [x]^{\perp_{\eta}\perp_{\varphi}}, \forall x\in\Pi. \]

  Put $\cH_4:=[\Pi]\cap\cH_5$. Since, obviously
  $f_{\varphi}([x])=f_{\varphi}|_{\Pi}([x])$ for $[x]\in\cH_4$,
  we have
  \begin{multline}
    \label{comm}
    f_{\varphi}([x])\in\cH_4 \Leftrightarrow
    f_\varphi([x])\in f_\varphi([x])^{\perp_{\eta}}
    \Leftrightarrow
    [x]^{\perp_{\varphi}\perp_{\eta}}\in [x]^{\perp_{\varphi}}
    \\
    \Leftrightarrow
    [x]^{\perp_{\eta}\perp_{\varphi}}\in [x]^{\perp_{\varphi}}
    \Leftrightarrow
    [x]\in [x]^{\perp_{\eta}}\Leftrightarrow [x]\in\cH_4.
  \end{multline}
  Thus,
  the collineation $g=f_{\varphi}|_{[\Pi]}$
  bijectively maps points of $\cH_4$ into points of $\cH_4$ and,
  consequently, $f_{\varphi}$ stabilizes $\cH_4$ as a set.
  In particular,
  $f_{\varphi}(\cH_4\setminus [R]^{\perp_{\eta}})=
  \cH_4\setminus f_{\varphi}([R]^{\perp_{\eta}})$.

Since $\cH_4$  is a hyperplane section of $\cH_5$   there are just two
  possible values for $|\cH_4|$, namely $(q^3+1)(q^2+1)$
  and $q^2(q^3+1)+1$, so
    \begin{equation}\label{numero 1}
     |\cH_4|\geq\min((q^3+1)(q^2+1), q^2(q^3+1)+1)=
    q^5+q^2+1.
   \end{equation}
  Likewise,
  since $\cH_4\cap [R]^{\perp_{\eta}}$ is
  a Hermitian variety in a projective space of dimension $2$
  and $f_{\varphi}$ is bijective when restricted to $\cH_4$
  we have
  \begin{multline}\label{numero 2}
    |\cH_4\cap f_{\varphi}([R]^{\perp_{\eta}})|= |\cH_4\cap [R]^{\perp_{\eta}}|\leq
    \\
    \leq\max( (q^3+1), (q^3+q^2+1), (q^2+1))= q^3+q^2+1.
\end{multline}

Note also that the case $[R]^{\perp_{\eta}}=[\Pi]$ and $[R]\in\cH_5$ can not happen because
 $[R]\not\in[\Pi].$ So, we need to consider the following three subcases:
 \begin{enumerate}[({a}1)]
\item
  \fbox{$[R]^{\perp_{\eta}}\neq[\Pi]$ and $[R]\not\in \cH_5$}.\\
  Since $[R]\not\in\cH_5$, $A\leq \frac{(q^4-1)(q^2-1)}{q-1}$.
Let $[x]\in\cH_4$. Since $f_{\varphi}(\cH_4)=\cH_4$, we have that
$[x]=f_{\varphi}([p])$ for some $[p]\in\cH_4$. In particular,
if $[p]\in\cH_4\setminus [R]^{\perp_{\eta}}$,
which is equivalent to say that
$[x]\in f_{\varphi}(\cH_4\setminus [R]^{\perp_{\eta}})=
\cH_4\setminus f_{\varphi}([R]^{\perp_{\eta}})$,
then the line
joining $[p]$ and $[R]$ meets $\cH_5$ in $(q+1)$
points and all of these points have the same image under $f_{\varphi}$.
So, $[x]$ admits either $q$ or $q+1$ preimages in
$\cH_5$ distinct from itself, according to $[x]\in[\langle p,R\rangle]$
or $[x]\not\in[\langle p,R\rangle]$.
By definition of $C$ and using Equation~(\ref{numero 1}) and Equation~(\ref{numero 2}), we get
  \begin{multline*}
    C\geq (q^2-1)q(|\cH_4|-|\cH_4\cap f_{\varphi}(R^{\perp_{\eta}})|)\geq\\
    \geq (q^2-1)q(q^5-q^3)=
    q^8-2q^6+q^4.
  \end{multline*}

Plugging this in Equation~\eqref{nuova pesi m=5}, we obtain
that for any $q$:
\begin{multline}
\wt(\varphi)=\frac{(q^3-q)}{(q^2+1)}\mu_5-\frac{(q^3-q)}{(q^4-1)}A +\frac{q}{(q^4-1)}C\geq\\
\frac{q^3-q}{q^2+1}(q^2+1)(q^5+1)-(q^3-q)(q+1)
+\frac{q}{(q^4-1)} (q^8-2q^6+q^4)=\\
(q^8-q^6)+ q^2\frac{q^2-1}{q^2+1}(q^3-q^2-1)>q^8-q^6.
\end{multline}

\item \fbox{$[R]^{\perp_{\eta}}\neq[\Pi]$ and $[R]\in \cH_5$}.\\
  In this case $A\leq\frac{(q^4-1)(q^2-1)}{q-1}+(q^2-1)$.
  As in Case a.1), consider a point $[x]\in\cH_4$.
  Then there is $[p]\in\cH_4$ such that $f_{\varphi}([p])=[x]$.
  The line joining a point
  $[p]\in \cH_4\setminus [R]^{\perp_{\eta}}$,
  with $[R]$ meets $\cH_5$ in $(q+1)$ points;
  every point of this
  line distinct from  $[R]$ is a preimage of $[x]=f_{\varphi}([p])$, so
  $[x]$
  admits at least $(q-1)$
  preimages in $\cH_5$ distinct from itself (they are exactly
  $q-1$ if $[x]\in [\langle p,R\rangle]$ and $q$ otherwise).
  Since $[p]\in  \cH_4\setminus [R]^{\perp_{\eta}}$ if
  and only if $[x]\in f_{\varphi}(\cH_4\setminus [R]^{\perp_{\eta}})=
  \cH_4\setminus f_{\varphi}([R]^{\perp_{\eta}})$,
  we have that any point $[x]\in  \cH_4\setminus f_{\varphi}([R]^{\perp_{\eta}})$
  admits at least $(q-1)$ preimages distinct from itself.

Take now $[s]\in \cH_4\cap [R]^{\perp_{\eta}}.$ Then $[s]$ is
collinear with $[R]$ and the line joining $[s]$ and $[R]$ is contained
in $\cH_5$. So, by the same argument as before, the preimages of a point $[x]=f_{\varphi}([s])\in\cH_4\cap f_{\varphi}([R]^{\perp_{\eta}})$
distinct from $[x]$ are at least $q^2-1$.
Since $f_{\varphi}$ is bijective when
restricted to
$\cH_4$, we have by Equations~(\ref{numero 1}) and Equations~(\ref{numero 2})
\begin{multline*} |\cH_4\setminus f_{\varphi}([R]^{\perp_{\eta}})|
  =|\cH_4\setminus [R]^{\perp_{\eta}}|
  = |\cH_4|-|\cH_4\cap [R]^{\perp_{\eta}}|\geq\\
  (q^5+q^2+1)-(q^3+q^2+1)=q^5-q^3.
\end{multline*}
On the other hand $|\cH_4\cap f_{\varphi}([R]^{\perp_{\eta}})|$ is
a (possibly degenerate or doubly-degenerate) Hermitian variety in a projective space of projective dimension $2$, so
\begin{multline*}
  |\cH_4\cap f_{\varphi}([R]^{\perp_{\eta}})|=|\cH_4\cap [R]^{\perp_{\eta}}|\geq\\
  \geq\min( (q^3+1), (q^3+q^2+1), (q^2+1))
  =(q^2+1).
\end{multline*}
Hence
  \begin{multline*}
    C
    \geq (q^2-1)\left((q-1)(|\cH_4\setminus [R]^{\perp_{\eta}}|
      )+(q^2-1)|\cH_4\cap [R]^{\perp_{\eta}}|\right)\geq\\
    (q^2-1)\left((q-1)(q^5-q^3)+(q^2-1)(q^2+1)\right)=\\
    (q^2-1)^2(q^4-q^3+q^2+1).
  \end{multline*}
  So,
\begin{multline}
\wt(\varphi)= \frac{(q^3-q)}{(q^2+1)}\mu_5-\frac{(q^3-q)}{(q^4-1)}A +\frac{q}{(q^4-1)}C\geq\\
\frac{q^3-q}{q^2+1}(q^2+1)(q^5+1) \\
- \frac{q^3-q}{q^4-1}( \frac{q^4-1}{q-1} (q^2-1) +(q^2-1))+\\
\frac{q}{q^4-1}(q^2-1)^2(q^4-q^3+q^2+1)=\\
(q^8-q^6)+ \frac{q^2(q^5-2q^4+q^2-q+1)}{q^2+1}> q^8-q^6.
\end{multline}

\item \fbox{$[R]^{\perp_{\eta}}=[\Pi]$ and $[R]\not\in \cH_5$}.\\
Since $[R]\not\in\cH_5$,  $A\leq\frac{(q^4-1)(q^2-1)}{q-1}$. Note that $\cH_4=[\Pi]\cap \cH_5$ is a non-degenerate Hermitian variety in a $3$-dimensional projective space, so $|\cH_4|=(q^3+1)(q^2+1).$
  The points of $[\Pi]\cap\cH_5$ not fixed by
  $f_{\varphi}$ admit exactly one preimage in $\cH_5$; so
  \[C\geq(|\cH_4|-|[\Fix(f_{\varphi})]|) (q^2-1)=(q^2-1)(q^5-q).\]
Thus
\begin{multline}
  \wt(\varphi)= \frac{(q^3-q)}{(q^2+1)}\mu_5-\frac{(q^3-q)}{(q^4-1)}A +\frac{q}{(q^4-1)}C\geq q^8-q^6.
\end{multline}
\end{enumerate}

 \item
  \fbox{$q^2+q+1\leq |[\Fix(f_{\varphi})]|\leq q^2+q+2$}. \\

Recall $\dim(\Rad(\varphi))=1$, so $\mathrm{Im}(f_{\varphi})$ is contained in a subspace
  of (vector) dimension  $4$ disjoint from $[R]:=[\Rad(\varphi)].$

  By Equation~\eqref{eAA}, we have
  \[
    A\leq (q^2-1)\Large(|[\Fix(f_{\varphi})]\cap\cH_5| +|[\Rad(\varphi)]\cap\cH_{5}|\Large).
  \]
  As each point of $f_{\varphi}(\cH_5)\cap\cH_5$ not fixed by $f_{\varphi}$
  has at least one preimage distinct from itself in $\cH_5$ we have
  \[
    C\geq (q^2-1)\Large(|f_{\varphi}(\cH_5)\cap\cH_5|-
    |[\Fix(f_{\varphi})]\cap\cH_5|\Large). \]

 Since $|[\Fix(f_{\varphi})]|\geq q^2+q+1$, there exist
  at least three points $[v_1]$, $[v_2]$, $[v_3]$ such that
  the vectors $v_i$ are linearly independent and
  \[
    f_{\varphi}(v_1)=\alpha_1 v_1,\quad
    f_{\varphi}(v_2)=\alpha_2 v_2,\quad
    f_{\varphi}(v_3)=\alpha_3 v_3. \]

By repeating the same argument as in Lemma~\ref{v4c1new} it is not difficult to see that $\alpha_1^{q+1}=\alpha_2^{q+1}=\alpha_3^{q+1}$.
Then applying Lemma~\ref{condizione (*)} rephrased for three independent vectors (instead of four), there exists a fixed point combination of $v_1,v_2,v_3$ with non-null coordinates.

Put $\Sigma:=[\langle v_1,v_2,v_3\rangle].$

The semilinear collineation $f_{\varphi}|_{\Sigma}$ is involutory,
  since it fixes a frame and, consequently, the linear
  collineation $(f_{\varphi}|_{\Sigma})^2$ is the identity.

Using now the same argument as in the beginning of
  point a), we see that $[x]^{\perp_{\varphi}\perp_{\eta}}=[x]^{\perp_{\eta}\perp_{\varphi}}$ for any $x\in \Sigma\cap\cH_5.$

  Put $\cH'=\cH_5\cap\Sigma.$ Hence,   (see Equation~\eqref{comm})
  $f_{\varphi}(\cH')=\cH'$.


Note also that  by Lemma~\ref{nuovo lemma} and Corollary~\ref{cofix}, $\Sigma\cap [\Fix(f_{\varphi})]$
  is a Baer subplane $\mathfrak G$ of $\Sigma$.

    In detail, if $|[\Fix(f_{\varphi})]|=q^2+q+1$, then we have
  $[\Fix(f_{\varphi})]\subseteq\Sigma$ and
  $[\Fix(f_{\varphi})]\cap\cH_5=[\Fix(f_{\varphi})]\cap\cH'$.

  If $|[\Fix(f_{\varphi})]|=q^2+q+2$,
  then by Corollary~\ref{cofix} there is a further point
  fixed by $f_{\varphi}$ not in $\Sigma$; so
  we have
  \begin{equation}
    \label{efix}
    |[\Fix(f_{\varphi})]\cap\cH_5|\leq |[\Fix(f_{\varphi})]\cap\cH'|+1.
  \end{equation}

  There are only three possibilities for $\cH'$:
  \begin{enumerate}[({b}1)]
  \item
    $\cH'$ is a non-degenerate Hermitian curve and $|\cH'|=q^3+1$,
  \item $\cH'$ is the union of $(q+1)$ lines through a
  point and $|\cH'|=q^2(q+1)+1$,
\item $ \cH'$ is a line  and $|\cH'|=q^2+1$.
\end{enumerate}
  We consider the three possibilities separately:
  \begin{enumerate}[({b}1)]
  \item
      If \fbox{$|\cH'|=q^3+1$}, then $\cH'$ is a non-degenerate
      Hermitian curve
      and the intersections between
      a non-degenerate Hermitian curve and a Baer subgeometry
      $\mathfrak G$ have size at most $2q+1$; see~\cite{23qt}.
    So
    \[ A\leq (q^2-1)(2q+1)+2(q^2-1),\qquad
      C\geq (q^2-1)(q^3-2q). \]
  It follows
        \begin{multline*}
      \wt(\varphi)\geq \frac{(q^3-q)}{(q^2+1)}\mu_5-\frac{(q^3-q)}{(q^4-1)}A +\frac{q}{(q^4-1)}C\geq\\
      (q^8-q^6)+q^3-q^2-4q+1+\frac{6q-1}{q^2+1}>q^8-q^6.
          \end{multline*}
        \item    If \fbox{$|\cH'|=q^3+q^2+1$}, then
          $\cH'$ is a cone $\Lambda$ of vertex a point $[v]$
          and basis a Baer subline of $(q+1)$ points.
          The overall number of fixed points on $\cH'$ is
          at most $q^2+q+1$.
          Put $\cH'':=\cH_5\cap f(\cH_5)\cap
          [\langle\Fix(f_{\varphi})\rangle]$.
          We now distinguish a few subcases.
          \begin{enumerate}[(b2.i)]
          \item
            $\cH''\not\subseteq\cH'$.
            In particular, $\cH''$ is not contained in the
            plane $\Sigma$.
            Since $\cH'\cup\cH''\subseteq \cH_5\cap f_{\varphi}(\cH_5)$,
            the projective dimension of
            $\langle\cH_5\cap f_{\varphi}(\cH_5)\rangle$ is $3$.
            Then $\cH_5\cap f_{\varphi}(\cH_5)$ is
            the intersection of two Hermitian
            surfaces in a projective space $\PG(3,q^2)$
            and contains the cone $\Lambda$.
            In particular, by Lemma~\ref{IntH4},
            $\cH_5\cap f_{\varphi}(\cH_5)$ must then have size
            at least $(q^2-1)^2$;
            thus we get
            \[ A\leq (q^2-1)(q^2+q+3), \]
            \[  C\geq (q^2-1)(q^4-2q^2+1-q^2-q-2), \]
            whence
            \begin{multline*}
              \wt(\varphi)\geq \frac{(q^3-q)}{(q^2+1)}\mu_5-\frac{(q^3-q)}{(q^4-1)}A +\frac{q}{(q^4-1)}C\geq\\
              q^8-q^6+(q^3-q^2)-6q+1+\frac{7q-1}{q^2+1}\geq
      q^8-q^6.
    \end{multline*}
  \item
    Suppose now that $\cH''\subseteq\cH'$.

    Clearly, $\cH''\subseteq\cH'\subseteq\Sigma$ by definition
    of $\cH'$ and
    $[R]\not\in\Sigma$
    as $\Sigma$ is spanned by fixed points of $f_{\varphi}$.
    \begin{itemize}
    \item
      If $[R]\not\in\cH_5$, then
      \[ A\leq (q^2-1)(q^2+q+1),\quad
        C\geq (q^2-1)(q^3-q). \]
      So,
      \[
        \wt(\varphi)\geq \frac{(q^3-q)}{(q^2+1)}\mu_5-\frac{(q^3-q)}{(q^4-1)}A +\frac{q}{(q^4-1)}C\geq
      q^8-q^6.
    \]
  \item
    Suppose now $[R]\in\cH_5$.
    Put $\tau:=|[\Fix(f_{\varphi})]\cap\cH'|\leq (q^2+q+1)$. Then,
    \[ A\leq (q^2-1)(\tau+1). \]
     The line joining $[R]$ with any fixed point $[x]$
     in $\cH_5\cap f(\cH_5)$ meets $\cH_5$ in at least $2$ points;
    so it is either contained in $\cH_5$ or it meets $\cH_5$ in
    $(q+1)$ points, including $[R]$ and $[x]$. So, on the line
    $[\langle R, x\rangle]$ there are at least $(q-1)$ points
    whose image is $[x]$ which are different from $[x]$ (and $[R]$).
    Obviously, none of these points is fixed.
    Also, observe that $[R]^{\perp_{\eta}}$ meets the plane of $\cH'$ in at least a
    line $\ell$ (otherwise the plane would be contained in $[R]^{\perp_{\eta}}$).
    Since $\cH'$ is the union of $q+1$ lines through a point,
    $|\ell\cap\cH'|\geq 1$. Let $[y]\in\ell\cap\cH'$.
    Then $f([y])$ admits at least $(q^2-1)$ preimages
    different from itself.
    So
    \[
      C\geq  (q^2-1)(
      (q^3+q^2+1-\tau)+(q-1)(\tau-1)+(q^2-1)).
    \]
    It follows that
    \begin{multline*}
      \wt(\varphi)\geq \frac{(q^3-q)}{(q^2+1)}\mu_5-\frac{(q^3-q)}{(q^4-1)}A +\frac{q}{(q^4-1)}C\geq\\
     (q^8-q^6)+(q^3+q^2-2)-\tau\left(q+\frac{1}{q^2+1}-1\right)+\frac{q+2}{q+1}.
    \end{multline*}
    Clearly, this bound is minimum when $\tau$ is maximum,
    that is $\tau=q^2+q+1$, so
\[
      \wt(\varphi)\geq (q^8-q^6)+q^2-2+\frac{2}{q^2+1}>q^8-q^6.
    \]
  \end{itemize}
\end{enumerate}
\item \fbox{$\cH'$ is a line}. \\
  \label{b3}
  By Theorem~\ref{new}
  the number of fixed points on a line is at most $q+1$ and there
  is at most one extra fixed point not in $\cH'$. So
  \[ A\leq (q^2-1)(q+3),\qquad C\geq0. \]
  So,
    \begin{multline*}
      \wt(\varphi)\geq \frac{(q^3-q)}{(q^2+1)}\mu_5-\frac{(q^3-q)}{(q^4-1)}A +\frac{q}{(q^4-1)}C\geq\\
      (q^8-q^6)+(q^3-q^2-4q+2)+\frac{6q-2}{q^2+1}>q^8-q^6.
    \end{multline*}
  \end{enumerate}
\item \fbox{$|[\Fix(f_{\varphi})]|\leq 2q+2$}.\\
   We have
  \[ A\leq (q^2-1)(2q+3),\qquad C\geq 0. \]
  So,
    \begin{multline*}
      \wt(\varphi)\geq \frac{(q^3-q)}{(q^2+1)}\mu_5-\frac{(q^3-q)}{(q^4-1)}A +\frac{q}{(q^4-1)}C\geq\\
      (q^8-q^6)+(q^3-2q^2-4q+4)+\frac{6q-4}{q^2+1}>q^8-q^6.
    \end{multline*}
    This completes the analysis.
  \end{enumerate}
\end{proof}
  \begin{remark}
    \label{r3}
    
  For $q=2$ and $m=5$, an exhaustive computer search shows that the
  characterization of Corollary~\ref{codd}
  does not hold, as there are $24948$
  codewords of minimum weight $192$; $19008$ of
  these are associated with bilinear forms with radical
  of dimension $1$ while the remaining $5940$ are  associated with forms
  with radical of dimension $3$. The
  forms with $3$-dimensional radical are as those described in Corollary~\ref{codd}.
  Incidentally, the full list of weights for this code is $0,192,216,224,232,256$.
\end{remark}
\begin{remark}
\label{r4}
  Consider the case $m=5$ and let $H$ and $S$ be respectively the matrix
  of the Hermitian form $\eta$ and of an alternating bilinear form $\varphi$.
\[ H:=\begin{pmatrix}
  1 & 0 & 0 & 0 & 0 \\
  0 & 0 & 1 & 0 & 0 \\
  0 & 1 & 0 & 0 & 0 \\
  0 & 0 & 0 & 0 & 1 \\
  0 & 0 & 0 & 1 & 0 \\
\end{pmatrix},\qquad
S:=\begin{pmatrix}
  0 & 0 & 0 & 0 & 0 \\
  0 & 0 & 1 & 0 & 0 \\
  0 & -1 & 0 & 0 & 0 \\
  0 & 0 & 0 & 0 & 1 \\
  0 & 0 & 0 & -1 & 0 \\
\end{pmatrix}. \]
Then $f_{\varphi}$ fixes exactly $q^3+q^2+q+1$ points and these points span the
polar hyperplane $[\Pi]=[\Rad(\varphi)]^{\perp_{\eta}}$ which has equation $x_1=0$.
A direct computation shows that the weight of the corresponding codeword is
$q^8-q^6$.  
\end{remark}


\section*{\Large Important Notice}
	\thispagestyle{empty}
	\setcounter{page}{0}
\large
\noindent
The paper
\begin{itemize}
\item I. Cardinali, L. Giuzzi, ``Line Hermitian Grassmann codes and
  their parameters'', \emph{Finite Fields Appl.} {\bfseries 51}:407--432
(2018)
\end{itemize}
contains a mistake in the proof of the case $m=5$ of the main theorem; this
does not affect the statement of the theorem itself, but changes the
corresponding characterization of the minimum weight codewords.

This document contains a copy  of the original version of the paper on arXiv
\url{https://arxiv.org/abs/1706.10255v2},
follwed by a \emph{corrigendum} note for this case.

We thank A. Cossidente for having informed us of the problem.\\
\\

\noindent The authors